\documentclass[12pt,a4paper]{amsart}

\usepackage{xypic}

\makeatletter                                            
\renewcommand{\@begintheorem}[2]{                        
\rm \trivlist \item [\hskip \labelsep {\bf #2\ \ #1.}]   
                                }                        
\makeatother                                             

\makeatletter

\DeclareFontFamily{U}{cyr}{}
\DeclareFontShape{U}{cyr}{m}{n}{
  <5> wncyr5 <6> wncyr6 <7> wncyr7 <8> wncyr8 <9> wncyr9 <10->
wncyr10}{}
\DeclareMathAlphabet{\mathcyr}{U}{cyr}{m}{n}

\input cyracc.def
\input epsf

\newcommand{\ts}{\vspace{\baselineskip}\noindent{\bf Proof.}$\;\;$}
\newcommand{\ZZ}{{\bf Z}}
\newcommand{\QQ}{{\bf Q}}

\newcommand{\CC}{{\bf C}}

\newcommand{\PP}{{\bf P}}

\newcommand{\lra}{\longrightarrow}

\newcommand{\calD}{{\mathcal D}}
\newcommand{\cM}{{\mathcal M}}
\newcommand{\cO}{{\mathcal O}}

\newcommand{\sig}{\sigma}

\newcommand{\rmod}{{\mbox{$\;{\rm mod}\;$}}}

\newcommand{\tX}{{\tilde{X}}}
\newcommand{\bY}{{\bar{Y}}}

\newcommand{\biota}{{\bar{\iota}}}
\def\blfootnote{\xdef\@thefnmark{}\@footnotetext}
\begin{document}
\title{Nikulin involutions on K3 surfaces}
\author{Bert van Geemen and Alessandra Sarti}
\address{Dipartimento di Matematica, Universit\`a di Milano,
  via Saldini 50, I-20133 Milano, Italia}
\address{Institut f\"ur  Mathematik, Universit\"at Mainz,
Staudingerweg 9, 55099 Mainz, Germany. {\it Current address}: Dipartimento di Matematica, Universit\`a di Milano,
  via Saldini 50, I-20133 Milano, Italia }

\email{geemen@mat.unimi.it}\email{sarti@mathematik.uni-mainz.de, sarti@mat.unimi.it}

\begin{abstract}
We study the maps induced on cohomology by a Nikulin (i.e.\ a
symplectic) involution on a K3 surface.
We parametrize
the eleven dimensional irreducible components of the
moduli space of algebraic K3 surfaces with a Nikulin involution
and we give examples of the general K3 surface in various components. 
We conclude with some remarks
on Morrison-Nikulin involutions, these are Nikulin involutions which
interchange two copies of $E_8(-1)$ in the N\'eron Severi group.

\end{abstract}

\maketitle

\blfootnote{The second author is supported by DFG Research Grant SA 1380/1-1.}

In his paper \cite{Nikulin} Nikulin started the study of finite groups of automorphisms on K3 surfaces, in particular
those leaving the holomorphic  two form invariant, these are called {\it symplectic}. He proves that  when the group $G$ is cyclic and acts symplectically, then $G\cong \ZZ/n\ZZ$, $1\leq n\leq 8$. Symplectic automorphisms of K3 surfaces of orders three, five and seven  
are investigated in the paper \cite{GS}. Here we consider the case of $G\cong \ZZ/2\ZZ$, generated by a symplectic involution $\iota$. Such involutions are called {\it Nikulin involutions} (cf.\cite[Definition 5.1]{morrison}). A Nikulin involution on the K3 surface $X$ has eight fixed points, hence the quotient $\bar{Y}=X/\iota$ has eight nodes, by blowing them up one obtains a K3 surface $Y$.\\
In the paper \cite{morrison} Morrison studies such involutions on algebraic K3 surfaces with Picard number $\rho\geq 17$ and in particular on those surfaces whose N\'eron Severi group contains two copies of $E_8(-1)$. These K3 surfaces always admit a Nikulin involution which interchanges the two copies of $E_8(-1)$. We call such involutions {\it Morrison-Nikulin involutions}.\\
The paper of Morrison motivated us to investigate Nikulin involutions in general. After a study of the maps on the cohomology induced by the quotient map, in the second section we show that an algebraic K3 surface with a Nikulin involution has $\rho\geq 9$ and that the N\'eron Severi group contains a primitive sublattice isomorphic with $E_8(-2)$. Moreover if $\rho=9$ (the minimal possible) then the following two propositions are the central results in the paper:

\bigskip

{\bf Proposition \ref{neronseveri}}. 
Let $X$ be a K3 surface with a Nikulin involution $\iota$
and assume that the N\'eron Severi group $NS(X)$ of $X$ has rank nine.
Let $L$ be a generator of $E_8(-2)^\perp\subset NS(X)$ with $L^2=2d>0$
and let
$$
\Lambda_{2d}\,:=\ZZ L\oplus E_8(-2)\quad(\subset NS(X)).
$$
Then we may assume that $L$ is ample and:
\begin{enumerate}
\item in case $L^2\equiv 2\;\text{mod}\,4$ we have $\Lambda_{2d}=NS(X)$;
\item in case $L^2\equiv 0\;\text{mod}\,4$ we have that either
$NS(X)\cong\Lambda_{2d}$ or $NS(X)\cong\Lambda_{\widetilde{2d}}$
where $\Lambda_{\widetilde{2d}}$ is the unique even lattice
containing $\Lambda_{2d}$ with
$\Lambda_{\widetilde{2d}}/\Lambda_{2d}\cong\ZZ/2\ZZ$ and such that $E_8(-2)$ is a primitive
sublattice of $\Lambda_{\widetilde{2d}}$.
\end{enumerate}

\bigskip

{\bf Proposition \ref{families}}. 
Let $\Gamma=\Lambda_{2d}$, $d\in\ZZ_{>0}$ or $\Gamma=\Lambda_{\widetilde{2d}}$,
$d\in 2\ZZ_{>0}$. Then there exists a K3 surface $X$ with a
Nikulin involution $\iota$ such that
$
NS(X)\cong \Gamma
$
and $(H^2(X,\ZZ)^\iota)^\perp\cong E_8(-2)$.

The coarse moduli space of $\Gamma$-polarized K3 surfaces has dimension 11
and will be denoted by $\cM_{2d}$ if $\Gamma=\Lambda_{2d}$ and by
$\cM_{\widetilde{2d}}$ if $\Gamma=\Lambda_{\widetilde{2d}}$.

\bigskip

 Thus we classified all the algebraic K3 surfaces with Picard number nine with a Nikulin involution. 
For the proofs we use lattice theory and the surjectivity of the period map for K3 surfaces. We also study the $\iota^*$-invariant line bundle $L$ on the general member of each family, for example in Proposition \ref{iota on bundles} we decompose the space $\PP H^0(X,L)^*$ into $\iota^*$-eigenspaces. This result is fundamental for the description of the $\iota$-equivariant map $X\longrightarrow \PP H^0(X,L)^*$. In section three we discuss various examples of the general K3 surface in these moduli spaces,  recovering well-known classical geometry in a few cases. We also describe the quotient surface $\bar{Y}$.

In the last section we give examples of K3 surfaces with an elliptic
fibration and a Nikulin involution which is induced by translation
by a section of order two in the Mordell-Weil group of the fibration.
Such a family has only ten moduli, and the minimal resolution of the quotient K3 surface $Y$ 
is again a
member of the same family. By using elliptic fibrations we also give an example of K3 
surfaces with a Morrison-Nikulin involution.
These surfaces with involution are parametrized by
three dimensional moduli spaces.
The Morrison-Nikulin involutions have interesting applications towards the
Hodge conjecture for products of K3 surfaces (cf.\ \cite{morrison},
\cite{galluzzi}). In section \ref{hodge} we briefly discuss possible
applications of the more general Nikulin involutions.

\section{General results on Nikulin Involutions}\label{general}
\subsection{Nikulin's uniqueness result}
A Nikulin involution $\iota$ of a K3 surface $X$ is an automorphism of order
two such that $\iota^*\omega=\omega$ for all $\omega\in H^{2,0}(X)$.
That is, $\iota$ preserves the holomorphic two form and thus it is a
symplectic involution. Nikulin, \cite[Theorem 4.7]{Nikulin}, proved
that any abelian group $G$ which
acts symplectically on a K3 surface, has a unique, up to isometry,
action on $H^2(X,\ZZ)$.

\subsection{Action on cohomology}\label{action}
D.\ Morrison (\cite[proof of Theorem 5.7]{morrison},)
observed that there exist
K3 surfaces with a Nikulin involution
which acts in the following way
on the second cohomology group:
$$
\iota^*:
H^2(X,\ZZ)\cong U^3\oplus E_8(-1)\oplus E_8(-1)\longrightarrow H^2(X,\ZZ),
\qquad
(u,x,y)\longmapsto (u,y,x).
$$
Thus for {\em any} K3 surface $X$ with
a Nikulin involution $\iota$ there is an isomorphism
$H^2(X,\ZZ)\cong U^3\oplus E_8(-1)\oplus E_8(-1)$ such that $\iota^*$ acts
as above.

Given a free $\ZZ$-module $M$ with an involution $g$,
there is an isomorphism
$$
(M,g)\cong M_1^s\oplus M_{-1}^t\oplus M_p^r,
$$
for unique integers $r,s,t$ (cf.\ \cite{Reiner}), where:
$$
M_1:=(\ZZ,\iota_1=1),\qquad M_{-1}:=(\ZZ,\iota_{-1}=-1),\qquad
M_p:=\left(\ZZ^2,\iota_p=\left(
\begin{array}{cc}
0&1\\1&0
\end{array}\right)
\right).
$$
Thus for a Nikulin involution acting on $H^2(X,\ZZ)$ the invariants
are $(s,t,r)=(6,0,8)$.

\subsection{The invariant lattice}\label{invariant lattice}
The invariant sublattice is:
$$
H^2(X,\ZZ)^\iota\,\cong
\{(u,x,x)\in U^3\oplus E_8(-1)\oplus E_8(-1)\,\}\cong\,
U^3\oplus E_8(-2).
$$
The anti-invariant lattice is the lattice
perpendicular to the invariant sublattice:
$$
(H^2(X,\ZZ)^\iota)^\perp\,\cong
\{(0,x,-x)\in U^3\oplus E_8(-1)\oplus E_8(-1)\,\}\cong\, E_8(-2).
$$
The sublattices $H^2(X,\ZZ)^\iota$ and $(H^2(X,\ZZ)^\iota)^\perp$ are
obviously primitive sublattices of $H^2(X,\ZZ)$.

\subsection{The standard diagram}\label{standard diagram}
The fixed point set of a Nikulin involution consists of exactly eight points
(\cite[section 5]{Nikulin}).
Let $\beta:\tilde{X}\rightarrow X$
be the blow-up of $X$ in the eight fixed points of $\iota$.
We denote by $\tilde{\iota}$ the involution on
$\tilde{X}$ induced by $\iota$.
Moreover, let $\bar{Y}=X/\iota$ be the eight-nodal quotient
of $X$, and let $Y=\tilde{X}/\tilde{\iota}$
be the minimal model of $\bar{Y}$, so $Y$ is a K3 surface.
This gives the `standard diagram':
$$
\begin{array}{ccc} X&\stackrel{\beta}{\longleftarrow}&\tilde{X}\\
\downarrow&&\downarrow \pi\\
\bar{Y}&{\longleftarrow}& Y.
\end{array}
$$
We denote by $E_i$, $i=1,\ldots,8$ the exceptional divisors in $\tilde{X}$
over the fixed points of $\iota$ in $X$, and by
$N_i=\pi(E_i)$ their images in $Y$, these are $(-2)$-curves.

\subsection{The Nikulin lattice}
The minimal primitive sublattice of $H^2(Y,\ZZ)$ containing the $N_i$ is
called the Nikulin lattice $N$ (cf.\ \cite[section 5]{morrison}).
As $N_i^2=-2$, $N_iN_j=0$ for $i\neq j$, the Nikulin lattice
contains the lattice $<-2>^8$. The lattice $N$ has rank eight and is spanned
by the $N_i$ and a class $\hat{N}$:
$$
N=\langle N_1,\ldots,N_8,\hat{N}\rangle,\qquad
\hat{N}:=(N_1+\ldots+N_8)/2.
$$
A set of $8$ rational curves on a K3 surface whose sum is divisible by $2$ in
the N\'eron Severi group is called an even set, see \cite{barth} and section \ref{examples}
for examples.

\subsection{The cohomology of $\tilde{X}$}
It is well-known that
$$
H^2(\tX, \ZZ)\cong H^2(X,\ZZ)\oplus (\oplus_{i=1}^8\ZZ E_i)
\cong U^3\oplus E_8(-1)^2\oplus  <-1>^8.
$$

For a smooth surface $S$ with torsion free $H^2(S,\ZZ)$,
the intersection pairing, given by the cup product to $H^4(S,\ZZ)=\ZZ$,
gives an isomorphism $H^2(S,\ZZ)\rightarrow Hom_\ZZ(H^2(S,\ZZ),\ZZ)$.

The map $\beta^*$ is:
$$
\beta^*:H^2(X,\ZZ)\lra H^2(\tX,\ZZ)=H^2(Y,\ZZ)\oplus (\oplus_{i=1}^8\ZZ E_i),
\qquad x\longmapsto(x,0),
$$
and its dual
$\beta_*:H^2(\tX,\ZZ)\rightarrow H^2(X,\ZZ)$ is $(x,e)\mapsto x$.

Let $\pi:\tX\rightarrow Y$ be the quotient map, let
$\pi^*:H^2(Y,\ZZ)\rightarrow H^2(\tX,\ZZ)$ be the induced map on the
cohomology and let
$\pi_*:H^2(\tX,\ZZ)\rightarrow H^2(Y,\ZZ)$ be its dual, so:
$$
\pi_*a\cdot b=a\cdot \pi^*b\qquad (a\in H^2(\tX,\ZZ),\;b\in H^2(Y,\ZZ)).
$$
Moreover, as $\pi^*$ is compatible with cup product we have:
$$
\pi^*b\cdot\pi^* c=2(b\cdot c)\qquad (b,c\in H^2(Y,\ZZ)).
$$

\subsection{Lattices}
For a lattice $M:=(M,b)$, where $b$ is a $\ZZ$-valued bilinear form on a free
$\ZZ$-module $M$, and an integer $n$ we let $M(n):=(M,nb)$. In particular,
$M$ and $M(n)$ have the same underlying $\ZZ$-module, but the identity map
$M\rightarrow M(n)$ is {\em not} an isometry unless $n=1$ or $M=0$.

\subsection{Proposition}\label{pi_*}
Using the notations and conventions as above,
the map $\pi_*:H^2(\tX,\ZZ)\longrightarrow H^2(Y,\ZZ)$
is given by
$$
\pi_*:U^3\oplus E_8(-1)\oplus E_8(-1)\oplus<-1>^8\longrightarrow
U(2)^3\oplus N \oplus E_8(-1)\;\hookrightarrow H^2(Y,\ZZ),
$$
$$
\pi_*:\;(u,x,y,z)\longmapsto (u,z, x+y).
$$
The map $\pi^*$, on the sublattice $U(2)^3\oplus N \oplus E_8(-1)$ of
$H^2(Y,\ZZ)$ is given by:
$$
\pi^*: U(2)^3\oplus N \oplus E_8(-1)\;\hookrightarrow\;
H^2(\tilde{X},\ZZ)\cong U^3\oplus E_8(-1)\oplus E_8(-1)\oplus<-1>^8,
$$
$$
\pi^*:(u,n,x)\longmapsto (2u,x,x,2\tilde{n}),
$$
here if $n=\sum n_iN_i$, $\tilde{n}=\sum n_iE_i$.

\ts
This follows easily from the results of Morrison.
In the proof of \cite[Theorem 5.7]{morrison},
it is shown that the image of each copy of $E_8(-1)$ under $\pi_*$
is isomorphic to $E_8(-1)$. As $E_8(-1)$ is unimodular, it is a
direct summand of the image of $\pi_*$. As $\pi_*\iota^*=\pi_*$,
we get that $\pi_*(0,x,0,0)=\pi_*(0,0,y,0)\in E_8(-1)$.
The $<-1>^8$ maps into $N$ (the image has index two).
As $U^3$ is a direct summand of $H^2(X,\ZZ)^{\iota}$, 
\cite[Proposition 3.2]{morrison} gives the first component.

As $\pi_*$ and $\pi^*$ are dual maps, $\pi^*a=b$ if for all
$c\in H^2(\tilde{X},\ZZ)$
one has $(b\cdot c)_{\tilde{X}}=(a\cdot \pi_*c)_Y$.
In particular, if $a\in U(2)^3$ and $c\in U^3$
we get $(\pi^*a\cdot c)_{\tilde{X}}=(a\cdot\pi_*c)_Y=2(a\cdot c)_{\tilde{X}}$
since we compute in $U(2)^3$,
hence $\pi^*a=2a$.
Similarly,
$(\pi^*N_i\cdot E_j)_{\tilde{X}}=(N_i\cdot \pi_*E_j)_Y=-2\delta_{ij}$,
so $\pi^*N_i=2E_i$ (this also follows from the fact that the $N_i$ are classes
of the branch curves,
so $\pi^*N_i$ is twice the class of $\pi^{-1}(N_i)=E_i$).
Finally for $x\in E_8(-1)$ and $(y,0)\in E_8(-1)^2$ we have
$(\pi^*x\cdot (y,0))_{\tilde{X}}=(x\cdot \pi_*(y,0))_Y=(x\cdot y)_Y$ and also
$(\pi^*x\cdot (0,y))_{\tilde{X}}=(x\cdot y)_Y$, so $\pi^*x=(x,x)\in E_8(-1)^2$.
\qed

\subsection{Extending $\pi^*$}
To determine the homomorphism $\pi^*:H^2(Y,\ZZ)\rightarrow H^2(\tX,\ZZ)$
on all of $H^2(Y,\ZZ)$, and not just on the sublattice of finite index
$U(2)^3\oplus N\oplus E_8(-1)$ we need to study the embedding
$U(2)^3\oplus N\hookrightarrow U^3\oplus E_8(-1)$.
This is done below.
For any
$x\in U^3\oplus E_8(-1)$, one has $2x\in U(2)^3\oplus N$
and $\pi^*(2x)$ determined as in Proposition \ref{pi_*}.
As $\pi^*$ is a homomorphism and lattices are torsion free, one finds
$\pi^*x$ as $\pi^*x=(\pi^*(2x))/2$.

\subsection{Lemma}\label{isotropic}
The sublattice of $(U(2)^3\oplus N)\otimes\QQ$ generated
by $U(2)^3\oplus N$ and the following six elements, each divided by two,
is isomorphic to $U^3\oplus E_8(-1)$:
$$
\begin{array}{lll}
e_1+(N_1+N_2+N_3+N_8),\quad &
                   e_2+(N_1+N_5+N_6+N_8),\quad &
                                           e_3+(N_2+N_6+N_7+N_8),\\
f_1+(N_1+N_2+N_4+N_8),\quad &
                   f_2+(N_1+N_5+N_7+N_8),\quad &
                                           f_3+(N_3+N_4+N_5+N_8),\\
\end{array}
$$
here $e_i,f_i$ are the standard basis of the $i$-th copy of $U(2)$ in $U(2)^3$.
Any embedding of $U(2)^3\oplus N$ into
$U^3\oplus E_8(-1)$ such that the image of $N$ is primitive in
$U^3\oplus E_8(-1)$ is isometric to this embedding.

\ts
The theory of embeddings of lattices can be found in \cite[section 1]{Nikulin bilinear}.
The dual lattice $M^*$ of a lattice $M=(M,b)$ is
$$
M^*=Hom(M,\ZZ)= \{x\in M\otimes\QQ:\; b(x,m)\in\ZZ\;\;\forall m\in M\}.
$$
Note that $M\hookrightarrow M^*$, intrinsically by $m\mapsto b(m,-)$
and concretely by $m\mapsto m\otimes 1$.
If $(M,b_M)$ and $(L,b_L)$ are lattices such that
$M\hookrightarrow L$, that is $b_M(m,m')=b_L(m,m')$ for $m,m'\in M$,
then we have a map $L\rightarrow M^*$ by $l\mapsto b_L(l,-)$.
In case $M$ has finite index in $L$, so
$M\otimes\QQ\cong L\otimes\QQ$, we get inclusions:
$$
M\hookrightarrow L\hookrightarrow L^*\hookrightarrow M^*.
$$
Therefore $L$ is determined by the image of $L/M$
in the finite group $A_M:=M^*/M$,
the discriminant group of $M$.

Since $b=b_M$ extends to a $\ZZ$-valued bilinear form on $L\subset M^*$
we get $q(l):=b_L(l,l)\in\ZZ$ for $l\in L$. If $L$ is an even lattice,
the discriminant form
$$
q_M:A_M\longrightarrow \QQ/2\ZZ,\qquad
m^*\longmapsto b_L(m^*,m^*)
$$
is identically zero on the subgroup $L/M\subset A_M$.
In this way one gets a bijection between even overlattices of
 $M$ and isotropic subgroups of $A_M$.
In our case $M=K\oplus N$, with $K=U(2)^3$,
so $A_M=A_K\oplus A_N$ and an isotropic subgroup of $A_M$ is the direct sum of an isotropic subgroup of $A_K$ and one isotropic subgroup of $A_N$.  We will see that $(A_K,q_K)\cong (A_N,-q_N)$, hence the even unimodular overlattices $L$ of $M$,
with $N$ primitive in $L$,
correspond to isomorphisms $\gamma:A_N\rightarrow A_K$ with
$q_N=-q_K\circ \gamma$. Then
one has that
$$
L/M=\{(\gamma(\bar{n}),\bar{n})\in A_M=A_K\oplus A_N:\; \bar{n}\in A_N\}.
$$
The overlattice $L_\gamma$ corresponding to $\gamma$ is:
$$
L_\gamma:=\{(u,n)\in K^* \oplus N^*:\;
\gamma(\bar{n})=\bar{u}\;\}.
$$
We will show that the isomorphism $\gamma$ is unique up to isometries of $K$ and $N$.\\
Let $e,f$ be the standard basis of $U$, so $e^2=f^2=0, ef=1$, then
$U(2)$ has the same basis with $e^2=f^2=0, ef=2$. Thus $U(2)^*$
has basis $e/2,f/2$ with $(e/2)^2=(f/2)^2=0, (e/2)(f/2)=2/4=1/2$.
Thus $A_K=(U(2)^*/U(2))^3\cong (\ZZ/2\ZZ)^6$, and the discriminant form $q_K$
on $A_K$ is given by
$$
q_K:A_K=(\ZZ/2\ZZ)^6\longrightarrow \ZZ/2\ZZ,\qquad
q_K(x)=x_1x_2+x_3x_4+x_5x_6.
$$

The Nikulin lattice $N$ contains $\oplus \ZZ N_i$ with $N_i^2=-2$,
hence $N^*\subset \ZZ (N_i/2)$. As $N=<N_i,(\sum N_i)/2>$
we find that $n^*\in\ZZ (N_i/2)$ is
in $ N^*$ iff $n^*\cdot (\sum N_i)/2\in\ZZ$, that is,
$n^*=\sum x_i(N_i/2)$ with $\sum x_i\equiv 0\rmod 2$. Thus we obtain an
identification:
$$
A_N=N^*/N=\{(x_1,\ldots,x_8)\in (\ZZ/2\ZZ)^8:\;\sum x_i=0\}/<(1,\ldots,1)>
\cong (\ZZ/2\ZZ)^6,
$$
where $(1,\ldots,1)$ is the image of $(\sum N_i)/2$.
Any element in $A_N$ has a unique representative which is either
$0$, $(N_i+N_j)/2$, with $i\neq j$
and $((N_i+N_j)/2)^2=1\rmod2\ZZ$,
or $(N_1+N_i+N_j+N_k)/2$ ($=(N_l+N_m+N_n+N_r)/2$), with distinct indices
and with $\{i,\ldots,r\}=\{2,\ldots,8\}$ and
$((N_1+N_i+N_j+N_k)/2)^2=0\rmod2$.
The quadratic spaces, over the field $\ZZ/2\ZZ$,
$((\ZZ/2\ZZ)^6,q_K)$ and
$((\ZZ/2\ZZ)^6,q_N)$ are isomorphic, an explicit
isomorphism is defined by
$$
\gamma:A_N\longrightarrow A_K,\qquad \gamma((N_1+N_2+N_3+N_8)/2)=e_1/2,
$$
etc.\ where we use the six elements listed in the lemma.

The orthogonal group of the quadratic space $((\ZZ/2\ZZ)^6,q_N)$
obviously contains $S_8$, induced by permutations of the basis
vectors in $(\ZZ/2\ZZ)^8$, and these groups are actually equal
cf. \cite{conway}.
Thus any two isomorphisms $A_N\rightarrow A_K$ preserving the quadratic
forms differ by an isometry of $A_N$ which is induced by a permutation of the
nodal classes $N_1,\ldots,N_8$.
A permutation of the $8$ nodal curves $N_i$ in $N$
obviously extends to an
isometry of $N$.

This shows that such an even unimodular overlattice of $U(2)^3\oplus N$
is essentially unique. As these are classified
by their rank and signature, the only possible one is $U^3\oplus E_8(-1)$.
Using the isomorphism $\gamma$, one obtains the lattice $L_\gamma$,
which is described in the lemma.
\qed

\subsection{The lattices $N\oplus N$ and $\Gamma_{16}$}\label{Gamma 16}
Using the methods of the proof of Lemma \ref{isotropic} we show that any
even unimodular overlattice $L$ of $N\oplus N$ such that $N\oplus\{0\}$ is
primitive in $L$, is isomorphic to $\Gamma_{16}(-1)$
(cf.\ \cite[Chapter V, 1.4.3]{Serre} ). The lattice $\Gamma_{16}(-1)$ is
the unique even unimodular negative definite lattice which is not generated
by its roots, i.e.\ by vectors $v$ with $v^2=-2$.

The discriminant form
$q_N$ of the lattice $N$ has values in $\ZZ/2\ZZ$, hence $q_N=-q_N$.
Therefore isomorphisms $\gamma:N\rightarrow N$ correspond to the
even unimodular overlattices $L_\gamma$ of $N\oplus N$ with
$N\oplus\{0\}$ primitive in $L_\gamma$. Since $N\oplus N$
is negative definite, so is $L_\gamma$. The uniqueness of the
overlattice follows,
as before, from the fact $O(q_N)\cong S_8$.
To see that this overlattice is $\Gamma_{16}(-1)$, recall that
$$
\Gamma_{16}=
\{x=(x_1,\ldots,x_{16})\in\QQ^{16}:
\;2x_i \in\ZZ,\;x_i-x_j\in\ZZ,\;\sum x_i\in 2\ZZ\;\},
$$
and the bilinear form on $\Gamma_{16}$ is given by $\sum x_iy_i$.
Let $e_i$ be the standard basis vectors of $\QQ^{16}$.
As
$$
N\oplus N\hookrightarrow\Gamma_{16}(-1),\qquad
(N_i,0)\longmapsto e_i+e_{i+8},\qquad (0,N_i)\longmapsto e_i-e_{i+8},
$$
is a primitive
embedding $N\oplus N$ into $\Gamma_{16}(-1)$ (note $(\hat{N},0)\mapsto
(\sum e_i)/2\in \Gamma_{16}$,  $(0,\hat{N})\mapsto
((\sum_{i=1}^8e_i)-(\sum_{i=9}^{16}e_i))/2\in\Gamma_{16}$)
the claim follows.

\section{Eleven dimensional families of K3 surfaces with a Nikulin involution}

\subsection{N\'eron Severi groups}\label {intro neronseveri}
As $X$ is a K3 surface it has $H^{1,0}(X)=0$ and 
$$
Pic(X)=NS(X)=H^{1,1}(X)\cap H^2(X,\ZZ)=
\{x\in H^2(X,\ZZ):\;x\cdot\omega=0\;\forall\omega\in H^{2,0}(X)\,\}.
$$
For $x\in (H^2(X,\ZZ)^\iota)^\perp$ we have $\iota^* x=-x$.
As $\iota^*\omega=\omega$ for $\omega\in H^{2,0}(X)$ we get:
$$
\omega\cdot x=\iota^*\omega\cdot \iota^* x=-\omega\cdot x\qquad
{\rm hence}\quad(H^2(X,\ZZ)^\iota)^\perp  \subset NS(X).
$$
As we assume $X$ to be algebraic, there is a very ample line bundle $M$
on $X$, so $M\in NS(X)$ and $M^2>0$. Therefore the N\'eron Severi group of
$X$  contains
$E_8(-2)\cong (H^2(X,\ZZ)^\iota)^\perp$ as a primitive sublattice
and has rank at least $9$.

The following proposition gives all even, rank $9$, lattices of signature
$(1+,8-)$ which contain $E_8(-2)$ as a primitive sublattice.
We will show in Proposition \ref{families} that any of these lattices
is the N\'eron Severi group of a K3 surface with a Nikulin involution.
Moreover, the moduli space of K3 surfaces, which contain such a
lattice in the N\'eron Severi group, is an $11$-dimensional complex variety.

\subsection{Proposition}\label{neronseveri}
Let $X$ be a K3 surface with a Nikulin involution $\iota$
and assume that the N\'eron Severi group of $X$ has rank $9$.
Let $L$ be a generator of $E_8(-2)^\perp\subset NS(X)$ with $L^2=2d>0$
and let
$$
\Lambda=\Lambda_{2d}\,:=\ZZ L\oplus E_8(-2)\quad(\subset NS(X)).
$$
Then we may assume that $L$ is ample and:
\begin{enumerate}
\item in case $L^2\equiv 2\;\text{mod}\,4$ we have $\Lambda=NS(X)$;
\item in case $L^2\equiv 0\;\text{mod}\,4$ we have that either
$NS(X)=\Lambda$ or $NS(X)\cong\tilde{\Lambda}$
where $\tilde{\Lambda}=\Lambda_{\widetilde{2d}}$ is the unique even lattice
containing $\Lambda$ with
$\tilde{\Lambda}/\Lambda\cong\ZZ/2\ZZ$ and such that $E_8(-2)$ is a primitive
sublattice of $\tilde{\Lambda}$.
\end{enumerate}

\ts
As $L^2>0$, either $L$ or $-L$ is effective, so may assume that
$L$ is effective. As there are no $(-2)$-curves in $L^\perp=E_8(-2)$,
any $(-2)$-curve $N$ has class $aL+e$ with $a\in\ZZ_{>0}$ and $e\in E_8(-2)$.
Thus $NL=aL^2>0$ and therefore $L$ is ample.

From the definition of $L$ and the description of the action of $\iota$ on
$H^2(X,\ZZ)$ it follows that $\ZZ L$ and $E_8(-2)$ respectively are primitive
sublattices of $NS(X)$. The discriminant group of $<L>$ is
$A_L:=<L>^*/<L>\cong \ZZ/2d\ZZ$ with generator $(1/2d)L$ where $L^2=2d$
and thus $q_L((1/2d)L)=1/2d$.
The discriminant group of $E_8(-2)$ is
$A_E\cong (1/2)E_8(-2)/E_8(-2)\cong(\ZZ/2\ZZ)^8$, as the quadratic form on
$E_8(-2)$ takes values in $4\ZZ$, the discriminant form $q_E$ takes values in
$\ZZ/2\ZZ$.

The even lattices $\tilde{\Lambda}$ which have
$\Lambda$ as sublattice of finite index correspond to isotropic subgroups
$H$ of $A_L\oplus A_E$ where $A_L:=<L>^*/<L>\cong \ZZ/2d\ZZ$. If $E_8(-2)$
is a primitive sublattice of $\tilde{\Lambda}$, $H$ must have
trivial intersection with both $A_L$ and $A_E$. Since $A_E$ is two-torsion,
it follows that $H$ is generated by $((1/2)L,v/2)$ for some $v\in E_8(-2)$.
As $((1/2)L)^2=d/2\;{\rm mod}\,2\ZZ$ and $(v/2)^2\in\ZZ/2\ZZ$,
for $H$ to be isotropic, $d$ must be even.
Moreover, if $d=4m+2$ we must have $v^2=8k+4$ for some $k$ and if
$d=4m$ we must have $v^2=8k$. Conversely, such a $v\in E_8(-2)$
defines an isotropic subgroup $<(L/2,v/2)>\subset A_L\oplus A_E$ which
corresponds to an overlattice $\tilde{\Lambda}$.
The group $O(E_8(-2))$ contains $W(E_8)$  (cf. \cite{conway}) which maps onto $O(q_E)$.
As $O(q_E)$ has three orbits on
$A_E$, they are $\{0\}$, $\{v/2:(v/2)^2\equiv 0\;(2)\}$
and $\{v/2:(v/2)^2\equiv 1\;(2)\}$, the overlattice is unique up to
isometry.
\qed

\subsection{Proposition}\label{families}
Let $\Gamma=\Lambda_{2d}$, $d\in\ZZ_{>0}$ or $\Gamma=\Lambda_{\widetilde{2d}}$,
$d\in 2\ZZ_{>0}$. Then there exists a K3 surface $X$ with a
Nikulin involution $\iota$ such that
$
NS(X)\cong \Gamma
$
and $(H^2(X,\ZZ)^\iota)^\perp\cong E_8(-2)$.

The coarse moduli space of $\Gamma$-polarized K3 surfaces has dimension 11
and will be denoted by $\cM_{2d}$ if $\Gamma=\Lambda_{2d}$ and by
$\cM_{\widetilde{2d}}$ if $\Gamma=\Lambda_{\widetilde{2d}}$.

\ts
We show that there exists a K3 surface $X$ with a Nikulin involution
$\iota$ such that $NS(X)\cong \Lambda_{\widetilde{2d}}$
and under this isomorphism
$(H^2(X,\ZZ)^\iota)^\perp\cong E_8(-2)$. The case $NS(X)\cong \Lambda_{2d}$
is similar but easier and is left to the reader.

The primitive embedding of
$\Lambda_{\widetilde{2d}}$ in the unimodular lattice
$U^3\oplus E_8(-1)^2$ is unique
up to isometry by \cite[Theorem 1.14.1]{Nikulin bilinear}, 
and we will identify $\Lambda_{\widetilde{2d}}$
with a primitive sublattice of $U^3\oplus E_8(-1)^2$ from now on.
We choose an
$\omega\in\Lambda_{\widetilde{2d}}^\perp\otimes_\ZZ\CC$ with
$\omega^2=0$, $\omega\bar{\omega}>0$ and general with these properties,
hence $\omega^\perp\cap (U^3\oplus E_8(-1)^2)=\Lambda_{\widetilde{2d}}$.
By the `surjectivity of the period map', there exists
a K3 surface $X$ with an isomorphism
$H^2(X,\ZZ)\cong U^3\oplus E_8(-1)^2$ such that
$NS(X)\cong\Lambda_{\widetilde{2d}}$.

The involution of $\Lambda=\ZZ L\oplus E_8(-2)$ which
is trivial on $L$ and $-1$ on $E_8(-2)$, extends to an involution of
$\Lambda_{\widetilde{2d}}=\Lambda+\ZZ(L/2,v/2)$.
The involution is trivial on the discriminant
group of $\Lambda_{\widetilde{2d}}$ which is isomorphic to $(\ZZ/2\ZZ)^6$.
Therefore it extends to an involution $\iota_0$
of $U^3\oplus E_8(-1)^2$ which is trivial on $\Lambda_{\widetilde{2d}}^\perp$.
As $((U^3\oplus E_8(-1)^2)^{\iota_0})^{\perp}=E_8(-2)$
is negative definite, contains no $(-2)$-classes and is contained in $NS(X)$,
results of Nikulin (\cite[Theorems 4.3, 4.7, 4.15]{Nikulin}) show
that $X$ has a Nikulin involution $\iota$ such that $\iota^*=\iota_0$
up to conjugation by an element of the Weyl group of $X$. Since we assume
$L$ to be ample and the ample cone is a fundamental domain for the
Weyl group action, we do get $\iota^*=\iota_0$, hence
$(H^2(X,\ZZ)^\iota)^\perp\cong E_8(-2)$.

For the precise definition of $\Gamma$-polarized K3 surfaces we refer to
\cite{dolgachev}. We just observe that each point of the moduli space corresponds
to a K3 surface $X$ with a primitive embedding
$\Gamma\hookrightarrow NS(X)$. The moduli space is a quotient of the
$11$-dimensional domain
$$
\calD_\Gamma=\{\omega\in \PP(\Gamma^\perp\otimes_\ZZ\CC):\;\omega^2=0,\;
\omega\bar{\omega}>0\,\}
$$
by an arithmetic subgroup of $O(\Gamma)$.
\qed

\subsection{Note on the Hodge conjecture}\label{hodge}
For a smooth projective surface $S$ with torsion free $H^2(S,\ZZ)$, let
$T_S:=NS(S)^\perp\;\subset H^2(S,\ZZ)$ and let
$T_{S,\QQ}=T_S\otimes_\ZZ\QQ$.
Then $T_S$, the transcendental lattice of $S$,
is an (integral, polarized) weight two Hodge structure.

The results in section \ref{general}
show that $\pi_*\circ \beta^*$ induces an isomorphism of rational
Hodge structures:
$$
\phi_\iota:T_{X,\QQ}\stackrel{\cong}{\lra} T_{Y,\QQ},
$$
in fact, both are isomorphic to $T_{\tilde{X},\QQ}$.
Any homomorphism of rational Hodge structures
$\phi:T_{X,\QQ} \rightarrow T_{Y,\QQ}$ defines,
using projection and inclusion,
a map of Hodge structures
$H^2(X,\QQ)\rightarrow T_{X,\QQ} \rightarrow T_{Y,\QQ}
\hookrightarrow H^2(Y,\QQ)$ and thus it gives
a Hodge (2,2)-class
$$
\phi\in H^2(X,\QQ)^*\otimes H^2(Y,\QQ)\cong H^2(X,\QQ)\otimes H^2(Y,\QQ)
\hookrightarrow H^4(X\times Y,\QQ),
$$
where we use Poincar\'e duality and the K\"unneth formula.
Obviously, the
isomorphism $\phi_\iota:T_{X,\QQ}\rightarrow T_{Y,\QQ}$
corresponds to the class of the codimension two cycle
which is the image of $\tX$ in $X\times Y$ under $(\beta,\pi)$.

Mukai showed that any homomorphism between $T_{S,\QQ}$ and $T_{Z,\QQ}$
where $S$ and $Z$ are K3 surfaces which is moreover an isometry
(w.r.t.\ the quadratic forms induced by the intersection forms)
is induced by an algebraic cycle if $\dim T_{S,\QQ}\leq 11$
(\cite[Corollary 1.10]{mukai}). Nikulin, \cite[Theorem 3]{nikulin cor},
strengthened this result and showed that it suffices that $NS(X)$ contains
a class $e$ with $e^2=0$. In particular, this implies that
any Hodge isometry $T_{S,\QQ}\rightarrow T_{Z,\QQ}$ is induced by an
algebraic cycle if $\dim T_{S,\QQ}\leq 18$ (cf.\ \cite[proof of Theorem 3]{nikulin cor}).

The Hodge conjecture predicts that any homomorphism of Hodge structures
between $T_{S,\QQ}$ and $T_{Z,\QQ}$ is induced by an algebraic cycle,
{\em without} requiring that it is an isometry.
There are few results in this direction, it is therefore maybe worth
noticing that $\phi_\iota$ is not an isometry if $T_X$ has odd rank,
see the proposition below. In \cite{galluzzi} a similar
result of D.\ Morrison in a more special case is used to obtain
new results on the
Hodge conjecture. In Proposition \ref{quotient fibration}
we show that there exists a K3 surface $X$ with Nikulin involution
where $T_{X,\QQ}$ has even rank and
$T_{X,\QQ}$ {\em is} isometric to $T_{Y,\QQ}$.

\subsection{Proposition} Let
$\phi_\iota:T_{X,\QQ}\stackrel{\cong}{\lra} T_{Y,\QQ}$ be the isomorphism
of Hodge structures induced by the Nikulin involution $\iota$ on $X$
and assume that $\dim T_{X,\QQ}$ is an odd integer.
Then $\phi_\iota$ is not an isometry.

\ts
Let $Q:\QQ^n\rightarrow \QQ$ be a quadratic form, then $Q$ is defined
by an $n\times n$ symmetric matrix, which we also denote by $Q$:
$Q(x):={}^txQx$. An isomorphism $A:\QQ^n\rightarrow \QQ^n$ gives an
isometry between $(\QQ^n,Q)$ and $(\QQ^n,Q')$ iff $Q'={}^tA^{-1}QA^{-1})$.
In particular, if $(\QQ^n,Q)\cong(\QQ^n,Q')$
the quotient $\det(Q)/\det(Q')$ must be a square in $\QQ^*$.

For a $\ZZ$-module $M$ we let $M_\QQ:=M\otimes_\ZZ\QQ$.
Let $V_X$ be the orthogonal complement of $E_8(-2)_\QQ\subset NS(X)_\QQ$,
then $\det(NS(X)_\QQ)=2^8\det(V_X)$ up to
squares. Let $V_Y$ be the orthogonal complement of $N_\QQ\subset NS(Y)_\QQ$
then $\det(NS(Y)_\QQ)=2^6\det(V_Y)$ up to
squares. Now $\beta_*\pi^*:H^2(Y,\QQ)\rightarrow H^2(X,\QQ)$
induces an isomorphism $V_X\rightarrow V_Y$ which
satisfies $(\beta_*\pi^* x)(\beta_*\pi^* y)=2xy$ for $x,y\in V_Y$,
hence $\det(V_X)=2^d\det(V_Y)$ where $d=\dim V_X=22-8-\dim T_{X,\QQ}$,
so $d$ is odd by assumption.

For a K3 surface $S$, $\det(T_{S,\QQ})=-\det(NS(S)_\QQ)$ and thus
$\det(T_{X,\QQ})/\det(T_{Y,\QQ})=2^{d+2}$ up to squares. As $d$ is odd and
$2$ is not a square in the multiplicative group of $\QQ$,
it follows that there exists no isometry between $T_{X,\QQ}$ and $T_{Y,\QQ}$.
\qed

\subsection{The bundle $L$}\label{bundle L}
In case $NS(X)$ has rank $9$, the ample generator
$L$ of $E_8(-2)^\perp\subset NS(X)$
defines a natural map
$$
\phi_L:X\longrightarrow \PP^g, \qquad g=h^0(L)-1=L^2/2+1
$$
which we will use to study $X$ and $Y$.
As $\iota^*L\cong L$, the involution $\iota$ acts as an involution
on $\PP^g=|L|^*$ and thus it has two fixed spaces $\PP^a,\PP^b$ with
$(a+1)+(b+1)=g+1$. The fixed points of $\iota$ map to these fixed spaces.
Even though $L$ is $\iota$-invariant, it is not the case in general that on
$\tX$ we have $\beta^*L=\pi^*M$
for some line bundle $M\in NS(Y)$. In fact, $\beta^*L=\pi^*M$ implies
$L^2=(\beta^*L)^2=(\pi^*M)^2=2M^2$ and as $M^2$ is even we get $L^2\in 4\ZZ$.
Thus if $L^2\not\in 4\ZZ$, the $\iota$-invariant line bundle $L$ cannot be
obtained by pull-back from $Y$.
On the other hand, if for example $|L|$ contains a reduced
$\iota$-invariant divisor $D$ which does not pass through the fixed points,
then $\beta^*D=\beta^{-1}D$ is invariant under $\tilde{\iota}$ on $\tX$
and does not contain any of the $E_i$ as a component.
Then $\beta^*D=\pi^*D'$ where $D'\subset Y$ is the reduced divisor
with support $\pi(\beta^{-1}D)$.

The following lemma collects the basic facts on $L$ and the splitting of
$\PP^g=\PP H^0(X,L)^*$.

\subsection{Proposition}\label{iota on bundles}
\begin{enumerate}
\item Assume that $NS(X)=\ZZ L\oplus E_8(-2)$.
Let $E_1,\ldots,E_8$ be the exceptional divisors on
$\tX$.

In case $L^2=4n+2$, there exist line bundles $M_1,M_2\in NS(Y)$
such that for a suitable numbering of these $E_i$ we have:
$$
\beta^*L-E_1-E_2=\pi^*M_1,\qquad
\beta^*L-E_3-\ldots-E_8=\pi^*M_2.
$$
The decomposition of $H^0(X,L)$ into $\iota^*$-eigenspaces is:
$$
H^0(X,L)\cong \pi^*H^0(Y,M_1)\oplus \pi^*H^0(Y,M_2),\qquad
(h^0(M_1)=n+2,\;h^0(M_2)=n+1).
$$
and the eigenspaces $\PP^{n+1}, \PP^{n}$
contain six, respectively two, fixed points.

In case $L^2=4n$, for a suitable numbering of the $E_i$ we have:
$$
\beta^*L-E_1-E_2-E_3-E_4=\pi^*M_1,\qquad
\beta^*L-E_5-E_6-E_7-E_8=\pi^*M_2
$$
with $M_1,M_2\in NS(Y)$.
The decomposition of $H^0(X,L)$ into $\iota^*$-eigenspaces is:
$$
H^0(X,L)\cong \pi^*H^0(Y,M_1)\oplus \pi^*H^0(Y,M_2),\qquad
(h^0(M_1)=h^0(M_2)=n+1).
$$
and each of the eigenspaces $\PP^n$ contains four fixed points.

\item
Assume that $\ZZ L\oplus E_8(-2)$ has index two in $NS(X)$.
Then there is a line bundle
$M\in NS(Y)$ such that:
$$
\beta^*L\cong\pi^*M,\qquad H^0(X,L)\cong H^0(Y,M)\oplus H^0(Y,M-\hat{N}),
$$
where $\hat{N}=(\sum_{i=1}^8 N_i)/2\in NS(Y)$ and this is the decomposition
of $H^0(X,L)$ into $\iota^*$-eigenspaces.
One has 
$h^0(M)=n+2$, $h^0(M-\hat{N})=n$,
and all fixed points
map to the eigenspace $\PP^{n+1}\subset \PP^{2n+1}=\PP^g$.
\end{enumerate}

\ts
The primitive embedding of
$\ZZ L \oplus E_8(-2)$ in the unimodular lattice
$U^3\oplus E_8(-1)^2$ is unique
up to isometry by \cite[Theorem 1.14.1]{Nikulin bilinear}.
Therefore if $L^2=2r$ we may assume that
$L=e_1+rf_1\in U\subset U^3\oplus E_8(-1)^3$
where $e_1,f_1$ are the standard basis of the first copy of $U$.

In case $r=2n+1$, it follows from Lemma \ref{isotropic} that
$(e_1+(2n+1)f_1+N_3+N_4)/2\in NS(Y)$. By Proposition \ref{pi_*},
$M_1:=(e_1+(2n+1)f_1+N_3+N_4)/2-N_3-N_4$
satisfies $\pi^*M_1=\beta^*L-E_3-E_4$.
Similarly, let $M_2=(e_1+(2n+1)f_1+N_3+N_4)/2-\hat{N}\in NS(Y)$, then
$\pi^*M_2=\beta^*L-(E_1+E_2+E_5+\ldots+E_8)$.

Any two sections $s,t\in H^0(X,L)$ lie in the same $\iota^*$-eigenspace
iff the rational function $f=s/t$ is $\iota$-invariant. Thus
$s,t\in \pi^*H^0(Y,M_i)$ are $\iota^*$-invariant,
hence each of these two spaces
is contained in an eigenspace of $\iota^*$ in $H^0(X,L)$.
If both are in the same eigenspace, then this eigenspace would have a section
with no zeroes in the $8$ fixed points of $\iota$. But a $\iota$-invariant
divisor on $X$ which doesn't pass through any fixed point is the pull back of
divisor on $Y$, which contradicts that $L^2$ is not a multiple of $4$.
Thus the $\pi^*H^0(Y,M_i)$ are in distinct eigenspaces. A dimension
count shows that $h^0(L)=h^0(M_1)+h^0(M_2)$, hence the $\pi^*H^0(Y,M_i)$
are the eigenspaces.

In case $r=2n$, again by Lemma \ref{isotropic} we have
$(e_1+N_1+N_2+N_3+N_8)/2\in NS(Y)$.
Let $M_1:=nf_1+(e_1+N_1+N_2+N_3+N_8)/2-(N_1+N_2+N_3+N_8)$
then $\pi^*M_1=\beta^*L-(E_1+E_2+E_3+E_8)$.
Put $M_2=M_1+\hat{N}-(N_4+N_5+N_6+N_7)$, then
$\pi^*M_2=\beta^*L-(E_4+E_5+E_6+E_7)$.
As above, the $\pi^*H^0(Y,M_i)$, $i=1,2$, are contained in distinct
eigenspaces and a dimension count again shows that $h^0(L)=h^0(M_1)+h^0(M_2)$.

If $\ZZ L\oplus E_8(-2)$ has index two in $NS(X)$, the (primitive)
embedding of $NS(X)$ into $U^3\oplus E_8(-1)$ is still unique up to
isometry. Let $L^2=4n$.
Choose an $\alpha\in E_8(-1)$ with $\alpha^2=-2$ if $n$ is odd,
and $\alpha^2=-4$ if $n$ is even.
Let $v=(0,\alpha,-\alpha)\in E_8(-2)\subset U^3\oplus E_8(-1)^2$ and let
$L=(2u,\alpha,\alpha)\in U^3\oplus E_8(-1)^2$ where $u=e_1+(n+1)/2f_1$
if $n$ is odd and $u=e_1+(n/2+1)f_1$ if $n$ is even.
Note that $L^2=4u^2+2\alpha^2=4n$ and that
$(L+v)/2=(u,\alpha,0)\in U^3\oplus E_8(-1)^2$. Thus we get a primitive
embedding of $NS(X)\hookrightarrow U^3\oplus E_8(-1)^2$ which extends
the standard one of $E_8(-2)\subset NS(X)$. Proposition \ref{pi_*} shows
that $\beta^*L=\pi^*M$ with
$M=(u,0,\alpha)\in U^3(2)\oplus N\oplus E_8(-1)\subset H^2(Y,\ZZ)$.
For the double cover $\pi:\tX\rightarrow Y$ branched along
$2\hat{N}=\sum N_i$
we have as usual:
$
\pi_*\cO_{\tilde{X}}=\cO_{Y}\oplus \cO_Y(-\hat{N})
$
hence, using the projection formula:
$$
H^0(\tX,\pi^*M)\cong H^0(Y,\pi_*(\pi^*M\otimes\cO_\tX)\cong
H^0(Y,M)\oplus H^0(Y,M-\hat{N}).
$$
Note that the sections in $\pi^*H^0(Y,M-\hat{N})$ vanish on all
the exceptional divisors, hence the fixed points of $\iota$ map to
a $\PP^{n+1}$.
\qed

\section{Examples}\label{examples}

\subsection{}
In Proposition \ref{families} we showed that K3 surfaces with a Nikulin
involution are parametrized by eleven dimensional moduli spaces
$\cM_{2d}$ and $\cM_{\widetilde{4e}}$ with $d,e\in\ZZ_{>0}$.
For some values of $d,e$ we will now work out the geometry of the
corresponding K3 surfaces.
We will also indicate how to verify that the moduli spaces are
indeed eleven dimensional.

\subsection{The case $\cM_{2}$}
Let $X$ be a K3 surface with Nikulin involution $\iota$
and $NS(X)\cong\ZZ L\oplus E_8(-2)$ with $L^2=2$ and $\iota^*L\cong L$
(cf.\ Proposition \ref{families}).
The map $\phi_L:X\rightarrow \PP^2$
is a double cover of $\PP^2$ branched over a sextic curve $C$,
which is smooth since there are no $(-2)$-curves in $L^\perp$.
The covering involution will be denoted by $i:X\rightarrow X$.
The fixed point locus of $i$ is isomorphic to $C$.

As $i^*$ is $+1$ on $\ZZ L$, $-1$ on $E_8(-2)$ and $-1$ on $T_X$,
whereas $\iota^*$ is $+1$ on $\ZZ L$, $-1$ on $E_8(-2)$ and $+1$ on $T_X$,
these two involutions commute. Thus $\iota$ induces an involution
$\bar{\iota}_{\PP^2}$ on $\PP^2$
(which is $\iota^*$ acting on $\PP H^0(X,L)^*$) and in suitable coordinates:
$$
\bar{\iota}_{\PP^2}:\;(x_0:x_1:x_2)\longmapsto (-x_0:x_1:x_2).
$$
We have a commutative diagram
$$
\begin{array}{rcrcl}
C&\hookrightarrow&X&\stackrel{\iota}{\longrightarrow}&X\\
\cong\downarrow &&\downarrow\phi&& \downarrow \phi\\
C&\hookrightarrow&\PP^2&\stackrel{\bar{\iota}_{\PP^2}}{\longrightarrow}&
\PP^2=X/i.
\end{array}
$$

The fixed points of $\bar{\iota}_{\PP^2}$ are:
$$
(\PP^2){\bar{\iota}_{\PP^2}}=l_0\cup \{p\},\qquad l_0:\;x_0=0,\quad p=(1:0:0).
$$
The line $l_0$ intersects
the curve $C$ in six points, which are the images of six fixed
points $x_3,\ldots,x_8$ of $\iota$ on $X$.
Thus the involution $\iota$ induces an involution on $C\subset X$ with
six fixed points.
The other two fixed points
$x_1,x_2$ of $\iota$ map to
the point $p$, so $i$ permutes these two fixed points of $\iota$.
In particular, these two points are not contained in $C$ so
$p\not\in C\,(\subset\PP^2)$, which will be important in the moduli count
below.
 The inverse image $C_2=\phi^{-1}(l_0)$ is a genus two curve
in the system $|L|$. Both $\iota$ and $i$ induce the hyperelliptic involution
on $C_2$. By doing then the quotient by $\iota$, since this has six fixed points on $C_2$ we obtain a rational curve $C_0$.

To describe the eight nodal surface $\bar{Y}=X/\iota$, we use the involution
$\bar{i}_{\bar{Y}}$ of $\bar{Y}$ which is induced by $i\in Aut(X)$.
Then we have:
$$
Q:=\bar{Y}/\bar{i}_\bY\;\cong\; X/<\iota,i>\;\cong\;\PP^2/\bar{\iota}_{\PP^2}.
$$
This leads to the following diagrams of double covers and fixed point
sets:

\vspace*{0.3cm}
\quad
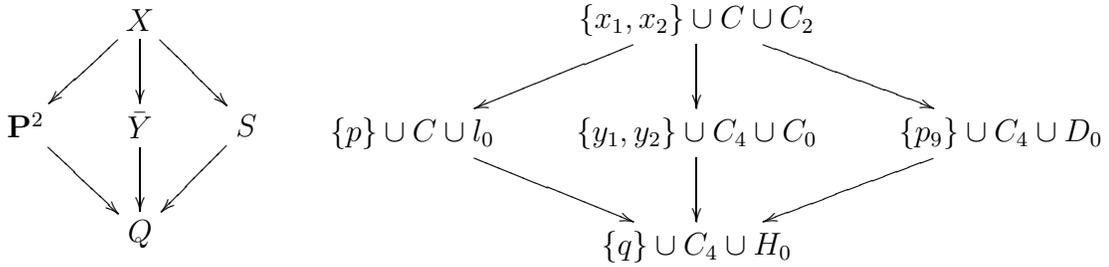
\begin{xy}
\xymatrix{
&X\ar[dl]\ar[d]\ar[rd]& \\
\PP^2\ar[rd]&\bar{Y}\ar[d]&S\ar[dl]\\
&Q&
}
\hspace*{5.0cm}
\xymatrix{
&\{x_1,x_2\}\cup C\cup C_2 \ar[dl]\ar[d]\ar[rd]& \\
\{p\}\cup C \cup l_0 \ar[rd]&
\{y_1,y_2\}\cup C_4 \cup C_0
\ar[d]&
\{p_9\}\cup C_4 \cup D_0
\ar[dl]\\
&\{q\}\cup C_4 \cup H_0&
}
\end{xy}
\vspace*{0.3cm}

The quotient of $\PP^2=X/i$ by $\bar{\iota}_{\PP^2}$ is isomorphic to a
quadric cone $Q$ in $\PP^3$ whose vertex $q$ is the image of the fixed point
$(1:0:0)$. In coordinates, the quotient map is:
$$
\PP^2\longrightarrow Q=\PP^2/\bar{\iota}_{\PP^2}\;\subset\PP^3,\qquad
(x_0:x_1:x_2)\longmapsto (y_0:\ldots:y_3)=(x_0^2:x_1^2:x_1x_2:x_2^2)
$$
and $Q$ is defined by $y_1y_3-y_2^2=0$.

The sextic curve $C\subset\PP^2$, which has genus $10$,
is mapped 2:1 to a curve $C_4$ on the cone. The double cover
$C\rightarrow C_4$ ramifies in the six points where
$C$ intersects the line $x_0=0$.
Thus the curve $C_4$ is smooth,
has genus four and degree six
(the plane sections of $C_4$ are the images of certain conic sections of
the branch sextic) and does not lie in a plane (so $C_4$ spans $\PP^3$).
The only divisor class $D$ of degree $2g-2$ with $h^0(D)\geq g$ on
a smooth curve of genus $g$ is the canonical class, hence
$C_4$ is a canonically embedded curve.
The image of the line $l_0$ is the plane section $H_0\subset Q$ defined
by $y_0=0$.

The branch locus in $Q$ of the double cover
$$
\bar{Y}\longrightarrow Q=\bar{Y}/\bar{i}_\bY
$$
is the union of two curves, $C_4$ and the plane section $H_0$,
these curves intersect in six points,
and the vertex $q$ of $Q$.

To complete the diagram, we consider the involution
$$
j:=\iota\circ i:X\longrightarrow X,\qquad S:=X/j.
$$
The fixed point set of $j$ is the (smooth) genus two curve $C_2$
lying over the line $l_0$ in $\PP^2$
(use $j(p)=p$ iff $\iota(p)=i(p)$ and consider the image of $p$ in $\PP^2$).
Thus the quotient surface $S$ is a smooth surface.
The Riemann-Hurwitz formula implies that the image of $C_2$ in $S$ is
a curve $D_0\in|-2K_{S}|$, note that $D_0\cong C_2$.

The double cover $S\rightarrow Q$ branches over the
curve $C_4\subset Q$ and the vertex $q\in Q$. It is well-known that
such a double cover is a Del Pezzo surface of degree 1
(\cite{demazure}, \cite{dolgort}) and the map $S\rightarrow Q\subset\PP^3$
is given by $\phi_{-2K}$,
which verifies that the image of $D_0$ is a plane section.

On the other hand, any
Del Pezzo surface of degree 1 is isomorphic to the blow up of $\PP^2$ in
eight points.
The linear system $|-K_S|$ corresponds to the pencil of elliptic
curves on the eight points, the ninth base point in $\PP^2$
corresponds to the unique
base point $p_9$ of $|-K_S|$ in $S$. The point $p_9$ maps to
the vertex $q\in Q$ under the 2:1 map
$\phi_{-2K}$ (\cite[p. 125]{dolgort}).
The N\'eron Severi group of $S$ is thus isomorphic to
$$
NS(S)\cong \ZZ e_0\oplus \ZZ e_1\oplus\ldots\oplus\ZZ e_8,\qquad
e_0^2=1,\quad e_i^2=-1\quad(1\leq i\leq 8)
$$
and $e_ie_j=0$ if $i\neq j$. The canonical class is
$K_S=-3e_0+e_1+\ldots +e_8$. Since $K_S^2=1$, we get a direct sum
decomposition:
$$
NS(S)\cong \ZZ K_S \oplus K_S^\perp\,\cong\, \ZZ K_S\,\oplus\,E_8(-1)
$$
(cf. \cite[VII.5]{dolgort}). The surface $S$ has 240 exceptional curves
(smooth rational curves $E$ with $E^2=-1$), cf.\ \cite[p.125]{dolgort}.
The adjunction formula
shows that $EK_S=-1$ and the map $E\mapsto E+K$ gives a bijection
between these exceptional curves and the roots of $E_8(-1)$, i.e.\ the
$x\in E_8(-1)$ with $x^2=-2$.
An exceptional divisor $E\subset S$ meets the branch curve $D_0\,(\in|-2K_S|)$
of $X\rightarrow S$ in two points, hence the inverse image of $E$ in $X$
is a $(-2)$-curve. Thus we get $240$ such $(-2)$-curves.
Actually,
$$
j^*:NS(S)=\ZZ K_S\,\oplus\,E_8(-1)\lra NS(X)=\ZZ L\oplus E_8(-2)
$$
is the identity on the $\ZZ$-modules and $NS(X)\cong NS(S)(2)$.
The class of such a $(-2)$-curve is
$L+x$, with $x\in L^\perp=E_8(-2)$, $x^2=-2$.
As $i^*(L+x)=L-x\neq L+x$, these $(-2)$-curves map pairwise
to conics in $\PP^2$, which must thus be tangent to the sextic $C$.
As also $\iota(L+x)=L-x$, these conics are invariant under $\biota_{\PP^2}$
and thus they correspond to plane sections of $Q\subset\PP^3$, tangent
to $C_4$, that is tritangent planes. This last incarnation of
exceptional curves in $S$ as tritangent
planes (or equivalently, odd theta characteristics of $C_4$) is of course
very classical.

Finally we compute the moduli. A $\bar\iota_{\PP^2}$-invariant
plane sextic which does not pass through $p=(1:0:0)$ has equation
$$
\sum a_{ijk}x_0^{2i}x_1^jx_2^k\qquad\qquad (2i+j+k=6,\; a_{000}\neq 0).
$$
The vector space spanned by such polynomials is 16-dimensional.
The subgroup of $GL(3)$ of elements commuting with $\biota_{\PP^2}$
(which thus preserve the eigenspaces) is isomorphic to $\CC^*\times GL(2)$,
hence the number of moduli is
$16-(1+4)=11$ as expected.

Alternatively, the genus four curves whose canonical image lies on a cone
have $9-1=8$ moduli (they have one vanishing even theta characteristic),
next one has to specify a plane in $\PP^3$, this gives again
$8+3=11$ moduli.

\subsection{The case $\cM_{6}$}
The map $\phi_L$ identifies $X$ with
a complete intersection of a cubic and a quadric in $\PP^4$.
According to Proposition \ref{iota on bundles}, in suitable coordinates
the Nikulin involution is induced by
$$
\iota_{\PP^4}:\PP^4\lra \PP^4,\qquad
(x_0:x_1:x_2:x_3:x_4)\longmapsto (-x_0:-x_1:x_2:x_3:x_4).
$$
The fixed locus in $\PP^4$ is:
$$
(\PP^4)^{\iota_{\PP^4}}=l\cup H,\qquad l:\;x_2=x_3=x_4=0,\quad H:\;x_0=x_1=0.
$$
The points $X\cap l$ and $X\cap H$ are fixed points of $\iota$ on $X$
and Proposition \ref{iota on bundles} shows that $\sharp(X\cap l)=2$,
$\sharp(X\cap H)=6$. In particular, the plane $H$ meets the quadric and cubic
defining $X$ in a conic and a cubic curve which intersect transversely.
Moreover, the quadric is unique, so must be invariant under $\iota_{\PP^4}$,
and, by considering the action of $\iota_{\PP^4}$ on the cubics
in the ideal of $X$,
we may assume that the cubic is invariant as well.
\begin{eqnarray*}
\begin{array}{lll}
l_{00}(x_2,x_3,x_4)x_0^2+l_{11}(x_2,x_3,x_4)x_1^2+l_{01}(x_2,x_3,x_4)x_0x_1
+f_3(x_2,x_3,x_4)&=&0\\
\alpha_{00}x_0^2+\alpha_{11}x_1^2+\alpha_{01}x_0x_1+f_2(x_2,x_3,x_4)&=&0\\
\end{array}
\end{eqnarray*}
where the $\alpha_{ij}$ are constants, the $l_{ij}$ are linear forms, and 
$f_2$, $f_3$ are
homogeneous polynomials of degree two and three respectively.
Note that the cubic contains the line $l:x_2=x_3=x_4=0$.

The projection from $\PP^4$ to the product of the eigenspaces 
$\PP^1\times\PP^2$ maps $X$ to a surface defined by an equation of bidegree 
$(2,3)$. In fact, the equations imply that $(\sum l_{ij}x_ix_j)/f_3=(\sum 
\alpha_{ij}x_ix_j)/f_2$ hence the image of $X$ is defined by the polynomial:
$(\sum l_{ij}x_ix_j)f_2-(\sum\alpha_{ij}x_ix_j)f_3$. Adjunction shows that a 
smooth surface of bidegree (2,3) is a K3 surface, so the equation defines 
$\bar{Y}$.
The space of invariant quadrics is $3+6=9$
dimensional and the space of cubics is $3\cdot 3+10=19$ dimensional.
Multiplying the quadric by a linear form $a_2x_2+a_3x_3+a_4x_4$
gives an invariant cubic.
The automorphisms of $\PP^4$ commuting with $\iota$ form a subgroup
which is isomorphic with $GL(2)\times GL(3)$ which has dimension $4+9=13$.
So the moduli space of such K3 surfaces has dimension:
\begin{eqnarray*}
(9-1)+(19-1)-3-(13-1)=11
\end{eqnarray*}
as expected.

\subsection{The case $\cM_4$}\label{degree four}
The map
$\phi_L:X\rightarrow \PP^3$ is an embedding whose image is
a smooth quartic surface.
From Proposition \ref{iota on bundles} the Nikulin involution
$\iota$ on $X\subset \PP^3\cong \PP(\CC^4)$ is induced by
$$
\tilde{\iota}:\CC^4\lra \CC^4,\qquad
(x_0,x_1,x_2,x_3)\longmapsto (-x_0,-x_1,x_2,x_3)
$$
for suitable coordinates.
The eight fixed points of the involution are the
points of intersection of these lines $x_0=x_1=0$ and $x_2=x_3=0$
with the quartic surface $X$.

A quartic surface which is invariant under $\tilde{\iota}$ and which does
not contain the lines has an equation which is
a sum of monomials $x_0^ax_1^bx_2^cx_3^d$ with $a+b=0,2,4$ and
$c+d=4-a-b$.

The quadratic polynomials invariant under $\tilde{\iota}$
define a map:
$$
\PP^3\lra \PP^5,\qquad (x_0:\ldots:x_3)\longmapsto
(z_0:z_1:\ldots:z_5)=(x_0^2:x_1^2:x_2^2:x_3^2:x_0x_1:x_2x_3)
$$
which factors over $\PP^3/\tilde{\iota}$.
Note that any quartic invariant monomial
is a monomial of degree two in the $z_i$. Thus if $f=0$ is the equation of
$X$, then $f(x_0,\ldots,x_3)=q(z_0,\ldots,z_5)$ for a quadratic form $q$.
This implies that
$$
\bar{Y}:\qquad q(z_0,\ldots,z_5)=0,\quad z_0z_1-z_4^2=0,\quad z_2z_3-z_5^2=0
$$
is the intersection of three quadrics.

The invariant quartics span a $5+9+5=19$-dimensional vector space.
On this space the subgroup $H$ of $GL(4)$ of elements which commute with
$\iota_{\PP^3}$ acts, it is easy to see that $H\cong GL(2)\times GL(2)$
(in block form).
Thus $\dim H=8$
and we get an $19-8=11$ dimensional family of quartic surfaces in
$\PP^3$, as desired. See \cite{Inose} for some interesting sub-families.

\subsection{The case $\cM_{\widetilde{4}}$}\label{falsodiv}
In this case $\ZZ L \oplus E_8(-2)$ has index two in $NS(X)$.
Choose a $v\in E_8(-2)$ with $v^2=-4$. Then we may assume that
$NS(X)$ is generated by $L,E_8(-2)$ and $E_1:=(L+v)/2$, cf.\ (the proof
of) Proposition \ref{neronseveri}.
Let $E_2:=(L-v)/2$, then $E_i^2=L^2/4+v^2/4=1-1=0$. By Riemann-Roch we have:
\begin{eqnarray*}
\chi(\pm E_i)=E_i^2/2+2=2
\end{eqnarray*}
and so $h^0(\pm E_i)\geq 2$ so $E_i$ or $-E_i$ is effective.
Now $L\cdot E_i=L^2/2+v/2 \cdot L=2$, hence $E_i$ is effective.
As $p_a(E_i)=1$ and $E_iN\geq 0$ for all $(-2)$-curves $N$,
each $E_i$ is the class of an elliptic fibration.
As $L=E_1+E_2$,
by \cite[Theorem 5.2]{saintdonat} the map $\phi_L$ is a 2:1 map
to a quadric $Q$ in $\PP^3$ and
it is ramified on a curve $B$ of bi-degree $(4,4)$.
The quadric is smooth, hence isomorphic to $\PP^1\times\PP^1$,
because there are no $(-2)$-curves in $NS(X)$ perpendicular to $L$.

Let $i:X\rightarrow X$ be the covering involution of
$X\rightarrow Q$. Then $i$ and the Nikulin-involution
$\iota$ commute.
The elliptic pencils $E_1$ and $E_2$ are permuted by
$\iota$ because $\iota^*L=L,\iota^*v=-v$.
This means that the involution $\bar{\iota}_Q$ on
$Q=\PP^1\times\PP^1$ induced
by $\iota$ acts as $((s:t),(u:v))\mapsto ((u:v),(s:t))$.
The quotient of $Q/\bar{\iota}_Q$ is well known to be isomorphic to
$\PP^2$.

The fixed point set of $\bar{\iota}_Q$ in $\PP^1\times\PP^1$ is
the diagonal $\Delta$. Thus $\Delta$ intersects the branch curve $B$
in eight points. The inverse image of these points in $X$ are the eight
fixed points of the Nikulin involution.

The diagonal maps to a conic $C_0$ in $\PP^2=Q/\bar{\iota}_Q$,
which gives the representation of a smooth quadric as double cover of
$\PP^2$ branched along a conic (in equations: $t^2=q(x,y,z)$).
The curve $B$ maps to a plane curve isomorphic to $\bar{B}=B/\iota$. As
$\iota$ has $8$ fixed points on the genus $9$ curve $B$,
the genus of $\bar{B}$ is $3$ and $\bar{B}\subset\PP^2$ is a quartic curve.

Let $j=i\iota=\iota i\in Aut(X)$. The fixed point set of $j$ is
easily seen to be the inverse image $C_3$ of the diagonal $\Delta\subset Q$.
As $C_3\rightarrow \Delta$ branches over the 8 points in $B\cap\Delta$,
$C_3$ is a smooth (hyperelliptic) genus three curve.
Thus the surface $S:=X/j$ is smooth and the image of $C_3$ in $S$
lies in the linear system $|-2K_S|$.
The double cover $S\rightarrow\PP^2$ is branched over the plane quartic
$\bar{B}\subset\PP^2$. This implies that $S$ is a Del Pezzo surface of degree
$2$, cf.\ \cite{demazure}, \cite{dolgort}.

This leads to the following diagrams of double covers and fixed point
sets:

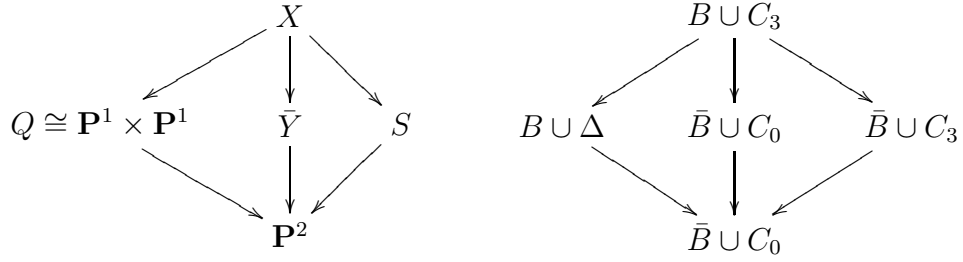
\vspace*{0.3cm}
\hspace*{2.0cm}
\begin{xy}
\xymatrix{
&X\ar[dl]\ar[d]\ar[rd]& \\
Q\cong\PP^1\times\PP^1\ar[rd]&\bar{Y}\ar[d]&S\ar[dl]\\
&\PP^2&
}
\hspace*{6.0cm}
\xymatrix{
& B \cup C_3 \ar[dl]\ar[d]\ar[rd]& \\
B \cup \Delta \ar[rd]&
\bar{B} \cup C_0
\ar[d]&
\bar{B} \cup C_3
\ar[dl]\\
&\bar{B} \cup C_0&
}
\end{xy}
\vspace*{0.3cm}

In particular the eight nodal surface $\bar{Y}$ is the double cover of
$\PP^2$ branched over the reducible sextic with components the conic $C_0$
and the quartic $\bar{B}$. The nodes of $\bar{Y}$ map to the intersection
points of $C_0$ and $\bar{B}$.

To count the moduli we note that the homogeneous polynomials of degree two
and four in three variables  span vector spaces of dimension $6$ and $15$,
as $\dim GL(3)=9$ we get:
$(6-1)+(15-1)-(9-1)=11$ moduli.

\subsection{The case $\cM_{8}$}\label{degree 8 new}
We have $H^0(X,L)\cong \pi^*H^0(Y,M_1)\oplus \pi^*H^0(Y,M_2)$
and $L^2=8$, $M_i^2=2$ so $h^0(L)=6$, $h^0(M_i)=3$ for $i=1,2$.
The image of $X$ under $\phi_L$ is the intersection of three quadrics
in $\PP^5$ and $\iota$ is induced by
$$
\tilde{\iota}:\CC^6\longrightarrow \CC^6,\qquad
(x_0,x_1,x_2,y_0,y_1,y_2)\longmapsto (x_0,x_1,x_2,-y_0,-y_1,-y_2).
$$
The multiplication map maps the
$21$-dimensional space $S^2H^0(X,L)$ onto the $18$-dimensional
space $H^0(X,2L)$. Using $\iota$ we can get some more information
on the kernel of this map, which are the quadrics defining $X\subset \PP^5$.
We have:
$$
S^2H^0(X,L) \cong
 \left(S^2H^0(Y,M_1)\oplus S^2H^0(Y,M_2)\right)\oplus
\left(H^0(Y,M_1)\otimes H^0(Y,M_2)\right),
$$
Moreover, as
$$
\beta^*(2L)=\pi^*M,\qquad{\rm with}\quad
M=2M_1+N_1+\ldots+N_4=2M_2+N_5+\ldots+N_8,
$$
(cf.\ Proposition \ref{iota on bundles})
we have the decomposition
$$
H^0(X,2L)\cong \pi^*H^0(Y,M)\oplus \pi^*H^0(Y,M-\hat{N}),\qquad
(h^0(M)=(M^2)/2+2=10,\;h^0(M-\hat{N})=8).
$$
In particular, the multiplication maps splits as:
$$
H^0(M_1)\otimes H^0(M_2)\lra H^0(Y,M-\hat{N})
$$
(vector spaces with dimensions with $3\cdot 3=9$ and $8$ resp.)
and
$$
S^2H^0(Y,M_1)\;\oplus\,S^2H^0(Y,M_2) \lra H^0(Y,M)
$$
(with dimensions $6+6=12$ and $10$ resp.).
Each of these two maps is surjective,
and as
$S^2H^0(Y,M_1)\rightarrow H^0(Y,M)$ is injective ($\phi_{M_1}$ maps
$Y$ onto $\PP^2$), the quadrics in the ideal of $X$ can be written as:
$$
Q_1(x)-Q_2(y)=0,\qquad Q_3(x)-Q_4(y)=0,\qquad B(x,y)=0
$$
with $Q_i$ homogeneous of degree two in three variables,
and $B$ of bidegree
$(1,1)$. Note that each eigenspace intersects $X$ in $2\cdot 2=4$ points.

The surface $\bar{Y}$ maps to $\PP^2\times \PP^2$ with the map
$\phi_{M_1}\times\phi_{M_2}$, its image is the image of $X$ under
the projections to the eigenspaces $\PP^5\rightarrow \PP^2\times\PP^2$.
As $(x_0:\ldots:y_2)\mapsto Q_1(x)/Q_2(y)$
is a constant rational function on $X$ and similarly
for $Q_3(x)/Q_4(y)$, there is a $c\in \CC$ such  that
the image of $X$ is contained in the
complete intersection of type $(2,2)$, $(1,1)$ in $\PP^2\times\PP^2$
defined by
$$
Q_1(x)Q_4(y)-cQ_3(x)Q_2(y)=0,\qquad B(x,y)=0.
$$
By adjunction, smooth complete intersections of this type are K3 surfaces.

To count the moduli, note that the first two equations come from a
$6+6=12$-dimensional vector space and the third comes from
a $3\cdot 3=9$-dimensional space.
The Grassmanian of $2$-dimensional subspaces of a $12$ dimensional space
has dimension $2(12-2)=20$. The subgroup of $GL(6)$ which commutes with
$\iota_{\PP^5}$ is isomorphic to $GL(3)\times GL(3)$ and has dimension
$9+9=18$. Thus we get
$20+(9-1)-(18-1)=11$ moduli, as expected.

\subsection{The case $\cM_{\widetilde{8}}$}\label{degree 8}
We have $H^0(X,L)\cong \pi^*H^0(Y,M)\oplus \pi^*H^0(Y,M-\hat{N})$
and $L^2=8$, $M^2=4$ so $h^0(M)=4,h^0(M-N)=2$.
The image of $X$ under $\phi_L$ is is the intersection of three quadrics
in $\PP^5$ and $\iota$ is induced by
$$
\tilde{\iota}:\CC^6\longrightarrow \CC^6,\qquad
(x_0,x_1,x_2,x_3,y_0,y_1)\longmapsto (x_0,x_1,x_2,x_3,-y_0,-y_1).
$$
To study the quadrics defining $X$, that is the kernel of the multiplication
map $S^2H^0(X,L)\rightarrow H^0(X,2L)$ we again split these spaces into
$\iota^*$-eigenspaces:
$$
S^2H^0(X,L)\cong
 \left(S^2H^0(Y,M)\oplus S^2H^0(Y,M-\hat{N})\right)\oplus
\left(H^0(Y,M)\otimes H^0(Y,M-\hat{N})\right),
$$
(with dimensions $21=(10+3)+8$)
and
$$
H^0(X,2L)\cong \pi^*H^0(Y,2M)\oplus \pi^*H^0(Y,2M-\hat{N})
$$
(with dimensions $h^0(2M)=10,\;h^0(2M-\hat{N})=8$).

This implies that there are no quadratic relations
in the $8$ dimensional space $H^0(Y,M)\otimes H^0(Y,M-\hat{N})$.
As $\phi_M$ maps $Y$ onto a quartic surface in $\PP^3$ and $M-\hat{N}$ is
a map of $Y$ onto $\PP^1$, the quadrics in the ideal of $X$ are of the form:
$$
y_0^2=Q_1(x),\qquad y_0y_1=Q_2(x),\qquad y_1^2=Q_3(x).
$$
The fixed points of the involution are the eight
points in the intersection of $X$ with the $\PP^3$ defined by $y_0=y_1=0$.

The image of $Y$ by $\phi_M$ is the image of the projection of $X$ from the
invariant line to the invariant $\PP^3$, which is
defined by $y_0=y_1=0$.
The image is the quartic surface
defined by $Q_1Q_3-Q_2^2=0$ which
can be identified with $\bar{Y}$. The equation is
the determinant of a symmetric $2\times 2$ matrix, which also implies
that this surface has $8$ nodes, (cf. \cite[Theorem 2.2]{Catanese}, \cite[section 3]{barth}),
the nodes form an even set (cf. \cite[Proposition 2.6]{Catanese}).

We compute the number of moduli.
Quadrics of this type span a space $U$ of dimension $3+10=13$.
The dimension of the Grassmanian of three dimensional subspaces of $U$ is
$3(13-3)=30$.
The group of automorphisms of $\CC^6$ which commute with $\iota_{\PP^5}$ is
$GL(2)\times GL(4)$.
So we have a $30-(4+16-1)=11$ dimensional space of such K3-surfaces
in $\PP^5$, as expected.

\subsection{The case $\cM_{12}$}\label{degree 12 new}
We have $H^0(X,L)\cong \pi^*H^0(Y,M_1)\oplus \pi^*H^0(Y,M_2)$
and $L^2=12$, $M_i^2=4$ so $h^0(L)=8$, $h^0(M_i)=4$ for $i=1,2$.
The image of $X$ under $\phi_L$ is the intersection of ten quadrics
in $\PP^7$. 

Following Example \ref{degree 8 new},
we use $\iota^*$ to split
the multiplication map from the
$36=(10+10)+16$-dimensional space $S^2H^0(X,L)$ onto the
$26=14+12$-dimensional
space $H^0(X,2L)$, again
$\beta^*(2L)=\pi^*M$ for an $M\in NS(Y)$ with $M^2=24$.
Thus we find $20-14=6$ quadrics of the type $Q_1(x)-Q_2(y)$
with $Q_i$ quadratic forms in $4$ variables, and $16-12=4$
quadratic forms $B_i(x,y)$, $i=1,\ldots,4$ where $x,y$ are coordinates
on the two eigenspaces in $H^0(X,L)$.

In particular, the projection from $\PP^7$ to the product of the eigenspaces
$\PP^3\times\PP^3$ maps $X$ onto a surface defined by $4$ equations of
bidegree $(1,1)$. Adjunction shows that a complete intersection of this
type is a K3 surface, so the
four $B_i$'s define $\bar{Y}\subset \PP^3\times\PP^3$.

Each $B_i$ can be written as: $B_i(x,y)=\sum_j l_{ij}(x)y_j$ with linear
forms $l_{ij}$ in $x=(x_0,\ldots,x_3)$.
The image of $\bar{Y}\subset\PP^3\times\PP^3$ under the projection
to the first factor is then defined by $\det(l_{ij}(x))=0$, which is a
quartic surface in $\PP^3$ as expected.
In fact, a point $x\in \PP^3$ has a non-trivial counter image
$(x,y)\in X\subset\PP^3\times\PP^3$ iff the matrix equation
$(l_{ij})y=0$ has a non-trivial solution.

As $X$ is not a complete intersection, we omit the moduli count.

\subsection{The case $\cM_{\widetilde{12}}$}\label{degree 12}
In this case $\beta^*L\cong \pi^*M$, $h^0(L)=8=5+3=h^0(M)+h^0(M-\hat{N})$.
We consider again the quadrics in the ideal of $X$ in Example \ref{degree 8}.
The space $S^2H^0(X,L)$ of quadrics on $\PP^7$ decomposes as:
$$
S^2H^0(X,L)\cong
 \left(S^2H^0(Y,M)+S^2H^0(Y,M-\hat{N})\right)\oplus
\left(H^0(Y,M)\otimes H^0(Y,M-\hat{N})\right),
$$
with dimensions $36=(15+6)+15$, whereas the sections of $2L$ decompose as:
$$
h^0(2L)=(4L^2)/2+2=26=14+12=
h^0(2M)\oplus h^0(2M-\hat{N}).
$$
Thus there are $(15+6)-14=7$ independent quadrics in the ideal
of $X\subset\PP^7$ which are
invariant and there are $15-12=3$
quadrics which are anti-invariant under the map
$$
\tilde{\iota}:\CC^8\lra \CC^8,\qquad
(x_0,\ldots,x_4,y_0,\ldots,y_2)\longmapsto (x_0,\ldots,x_4,-y_0,\ldots,-y_2).
$$
An invariant quadratic polynomial looks like
$q_0(x_0,\ldots,x_4)+q_1(y_0,y_1,y_2)$,
and since the space of quadrics in three
variables is only $6$ dimensional, there is one non-zero quadric
$q$ in the ideal of
the form $q=q(x_0,\ldots,x_4)$.
An anti-invariant quadratic polynomial is of bidegree
$(1,1)$ in $x$ and $y$.
In particular, the image of the projection of $X$ to the product of the
eigenspaces $\PP^4\times\PP^2$
is contained in one hypersurface
of bidegree $(2,0)$ and in three hypersurfaces of bidegree $(1,1)$.
The complete intersection of four general such hypersurfaces is a K3 surface
(use adjunction and $(2+3\cdot 1,3\cdot 1)=(5,3)$).

The three anti-invariant quadratic forms can be written as
$\sum_jl_{ij}(x)y_j$, $i=1,2,3$.
The determinant of the $3\times 3$ matrix of linear forms $(l_{ij}(x))$,
defines a cubic form which is an equation for the image of $X$ in $\PP^4$
(cf.\ Example \ref{degree 12 new}).
Thus the projection $\bar{Y}$ of $X$ to $\PP^4$ is the
intersection of the quadric defined by $q(x)=0$ and a cubic.

The projection to $\PP^2$ is $2$:$1$, as it should be, since for general
$y\in\PP^2$ the three linear forms in $x$ given by $\sum_j l_{ij}(x)y_j$
define a line in $\PP^4$ which cuts the quadric $q(x)=0$ in two points.

\section{Elliptic fibrations with a section of order two}

\subsection{Elliptic fibrations and Nikulin involutions}
\label{degree zero}\label{elliptic}
Let $X$ be a K3 surface which has an elliptic fibration
$f:X\rightarrow \PP^1$
with a section $\sig$. The set of sections of $f$ is a group,
the Mordell-Weil group $MW_f$, with identity element $\sig$.
This group acts on $X$ by translations and these translations preserve
the holomorphic two form on $X$.
In particular, if there is an element $\tau\in MW_f$ of
order two, then translation by $\tau$ defines a Nikulin involution $\iota$.

In that case the Weierstrass equation of $X$ can be put in the form:
$$
X:\qquad y^2=x(x^2+a(t)x+b(t)) 
$$
the sections $\sig,\tau$ are given by the section at infinity and
$\tau(t)=(x(t),y(t))=(0,0)$. For the general fibration on a K3 surface $X$,
the degrees of $a$ and $b$ are $4$ and $8$ respectively.

\subsection{Proposition}\label{quotient fibration}
Let $X\rightarrow \PP^1$ be a general elliptic fibration with sections
$\sigma,\tau$ as above in section \ref{elliptic}. and let $\iota$
be the corresponding Nikulin involution on $X$. These fibrations form a
$10$-dimensional family.

The quotient K3 surface $Y$ also has an elliptic fibration:
$$
Y:\qquad y^2=x(x^2-2a(t)x+(a(t)^2-4b(t)), 
$$

We have:
$$
NS(X)\cong NS(Y)\cong U\oplus N,\qquad T_X\cong T_Y\cong U^2\oplus N.
$$

The bad fibers of $X\rightarrow \PP^1$ are eight fibers of type $I_1$
(which are rational curves wit a node) over the zeroes of $a^2-4b$
and eight fibers of type $I_2$
(these fibers are the union of two $\PP^1$'s meeting in two points)
over the zeroes of $b$. The bad fibers of $Y\rightarrow \PP^1$ are
eight fibers of type $I_2$
over the zeroes of $a^2-4b$
and eight fibers of type $I_1$ over the zeroes of $b$.

\ts
Since $X$ has an elliptic fibration with a section, $NS(X)$ contains a copy
of the hyperbolic plane $U$
(with standard basis the class of a fiber $f$ and $f+\sigma$).
The discriminant of the Weierstrass model of $X$ is $\Delta_X=b^2(a^2-4b)$
and the fibers of the Weierstrass model over the zeroes of $\Delta_X$
are nodal curves.
Thus $f:X\rightarrow \PP^1$ has eight fibers of type $I_1$
(which are rational curves with a node) over the zeroes of $a^2-4b$
and 8 fibers of type $I_2$
(these fibers are the union of two $\PP^1$'s meeting in two points)
over the zeroes of $b$.

The components of the singular fibers which do not meet
the zero section $\sigma$, give a sublattice $<-2>^8$ perpendicular to $U$.
If there are no sections of infinite order, the lattice $U\oplus<-2>^8$
has finite index in the N\'eron Severi group of $X$. Hence $X$ has
$22-2-10=10$ moduli.
One can also appeal to \cite{Shimada} where the N\'eron Severi group of the
general elliptic K3 fibration with a section of order two is determined.
To find the moduli from the Weierstrass model, note that
$a$ and $b$ depend on $5+9=14$ parameters.
Using transformations of the type $(x,y)\mapsto (\lambda^2 x,\lambda^3 y)$
(and dividing the equation by $\lambda^6$)
and the automorphism group $PPG2)$ of $\PP^1$ we get $14-1-3=10$ moduli.

The Shioda-Tate formula (cf.\ e.g.\ \cite[Corollary 1.7]{shioda}) shows
that the discriminant of the N\'eron Severi
group is $2^8/n^2$ where $n$ is the order of the torsion subgroup of $MW_f$.
The curve defined by $x^2+a(t)x+b(t)=0$ cuts out the remaining pair of points
of order two on each smooth fiber. As it is irreducible
in general,
$MW_f$ must be cyclic.
If there were a section $\sig$ of order
four, it would have to satisfy $2\sig=\tau$.
But in a fiber of type $I_2$ the
complement of the singular points is the group $G=\CC^*\times(\ZZ/2\ZZ)$ and
the specialization $MW_f\rightarrow G$ is an injective homomorphism.
Now $\tau$ specializes to $(\pm 1,\bar{1})$ (the sign doesn't matter)
since $\tau$ specializes to the node in the Weierstrass model.
But there is no $g\in G$ with $2g=(\pm 1,\bar{1})$. We conclude that for
general $X$ we have $MW_f=\{\sigma,\tau\}\cong\ZZ/2\ZZ$
and that the discriminant
of the N\'eron Severi group of $X$ is $2^6$.

The N\'eron Severi group has $\QQ$ basis $\sigma,f,N_1,\ldots,N_8$
where the $N_i$ are the components of the $I_2$ fibers not meeting $\sigma$.
As $\tau\cdot \sigma=0$, $\tau\cdot f=1$ and $\tau\cdot N_i=1$, we get:
$$
\tau=\sigma+2f-\hat{N},\qquad \hat{N}=(N_1+\ldots+N_8)/2.
$$
Thus the smallest primitive sublattice containing the $N_i$ is the Nikulin
lattice. Comparing discriminants we conclude that:
$$
NS(X)=\langle s,f\rangle \oplus \langle N_1,\ldots,N_8,\hat{N}\rangle
\;\cong\; U\oplus N.
$$

The transcendental lattice $T_X$ of $X$ can be determined as follows.
It is a lattice of signature $(2+,10-)$ and its discriminant form
is the opposite of the one of $N$,
but note that $q_N=-q_N$ since $q_N$ takes values in $\ZZ/2\ZZ$.
Moreover, $T_X^*/T_X\cong N^*/N\cong (\ZZ/2\ZZ)^6$.
Using \cite[Corollary 1.13.3]{Nikulin bilinear}, we find that
$T_X$ is uniquely determined by the signature and the discriminant form.
The lattice $U^2\oplus N$ has these invariants, so
$$
T_X\cong U^2\oplus N.
$$

As the Nikulin involution preserves the fibers of the elliptic
fibration on $X$,
the desingularisation $Y$ of the quotient $X/\iota$ has
an elliptic fibration $g:Y\rightarrow \PP^1$,
with a section $\bar{\sig}$, (the image of $\sig$).
The Weierstrass equation of $Y$ can be found from \cite[p.79]{ST}.

The discriminant of the Weierstrass model of $Y$
is $\Delta_Y=4b(a^2-4b)^2$ and, reasoning as before, we find the bad fibers
of $g:Y\rightarrow \PP^1$. In particular,
the $I_1$ and $I_2$ fibers of $X$ and $Y$ are indeed `interchanged'.

Geometrically, the reason for this is as follows.
The fixed points of translation by $\tau$ are the eight nodes in the
$I_1$-fibers,
blowing them up gives $I_2$-type fibers which map to $I_2$-type fibers in $Y$.
The exceptional curves lie in the ramification locus of the quotient map,
the other components, which meet $\sigma$, map 2:1 to components
of the $I_2$-fibers which meet $\bar{\sigma}$.
The two components of an $I_2$-fiber in $X$ are interchanged and also the
two singular
points of the fiber are permuted, so in the quotient this gives an
$I_1$-type fiber.
\qed

\subsection{Remark} Note that $NS(X)\oplus T_X\cong U^3\oplus N^2$,
however, there is {\emph{no}} embedding of $N^2$ into $E_8(-1)^2$,
such that $N\oplus\{0\}\,(\subset NS(X))$ is primitive in $E_8(-1)^2$.
However, $N^2\subset \Gamma_{16}(-1)$ (cf.\ section \ref{Gamma 16}),
an even, negative definite, unimodular lattice of rank 16 and
$U^3\oplus \Gamma_{16}(-1)\cong U^3\oplus E_8(-1)^2$ by the classification
of even indefinite unimodular quadratic forms.

\subsection{Morrison-Nikulin involutions}
D.\ Morrison observed that a K3 surface $X$
having two perpendicular copies of
$E_8(-1)$ in the N\'eron Severi group has a Nikulin involution
which exchanges the two copies of $E_8(-1)$,
cf.\ \cite[Theorem 5.7]{morrison}.
We will call such an involution a Morrison-Nikulin involution.
This involution then has the further property that $T_Y\cong T_X(2)$
where $Y$ is the quotient K3 surface and we have a Shioda-Inose
structure on $Y$ (cf.\ \cite[Theorem 6.3]{morrison})

\subsection{Moduli}\label{moduli MN}
As $E_8(-1)$ has rank eight and is negative definite, a projective K3 surface
with a Morrison-Nikulin involution has a N\'eron Severi group of rank at least
$17$ and hence has at most three moduli. In case the N\'eron Severi group has
rank exactly $17$, we get
$$
NS(X)\cong \langle 2n\rangle \oplus E_8(-1)\oplus E_8(-1)
$$
since the sublattice $E_8(-1)^2$ is unimodular.
Results of Kneser and Nikulin, \cite[Corollary 1.13.3]{Nikulin bilinear},
guarantee that the transcendental lattice $T_X:=NS(X)^\perp$
is uniquely determined by its signature and discriminant form.
As the discriminant form of $T_X$ is the opposite of the one on $NS(X)$
we get
$$
T_X\cong \langle -2n\rangle \oplus U^2.
$$

In case $n=1$ such a three dimensional family can be obtained
from the double covers of $\PP^2$ branched along a sextic curve
with two singularities which are locally isomorphic to $y^3=x^5$.
The double cover then has two singular points of type $E_8$,
that is, each of these can be resolved by
eight rational curves with incidence graph $E_8$.
As the explicit computations are somewhat lengthy and involved,
we omit the details. See \cite{persson} and \cite{degtyarev}
for more on double covers of $\PP^2$ along singular sextics.

\subsection{Morrison-Nikulin involutions on elliptic fibrations}
We consider a family of K3 surfaces with an elliptic fibration
with a Morrison-Nikulin involution induced by translation by a section
of order two. It corresponds to the family with $n=2$ from section
\ref{moduli MN}.

Note that in the proposition below we describe a
$K3$ surface $Y$ with a Nikulin involution and quotient K3
surface $X$ such  that $T_Y=T_X(2)$, which is the `opposite' of
what would happen if the involution of $Y$ was a Morrison-Nikulin involution.
it is not hard to see that there is no primitive embedding
$T_Y\hookrightarrow U^3$, so $Y$ does not have a Morrison-Nikulin involution
at all (cf. \cite[Theorem 6.3]{morrison}).

\subsection{Proposition}\label{elliptic n-m}
Let $X\rightarrow \PP^1$ be a general elliptic fibration defined by the
Weierstrass equation
$$
X:\qquad y^2=x(x^2+a(t)x+1),\qquad a(t)=a_0+a_1t+a_2t^2+t^4
\in\CC[t].
$$
The K3 surface $X$ has a Morrison-Nikulin involution
defined by translation by the section, of order two,
$t\mapsto (x(t),y(t))=(0,0)$.
Then:
$$
NS(X)=\langle 4 \rangle\oplus E_8(-1)\oplus E_8(-1),\qquad
T_X=\langle -4\rangle\oplus U^2.
$$
The bad fibers of the fibration are nodal cubics
(type $I_1$) over the eight zeroes of $a^2(t)-4$ and one fiber
 of type $I_{16}$ over $t=\infty$.

The quotient K3 surface $Y$ has an elliptic fibration defined by the
Weierstrass model:
$$
Y:\qquad y^2=x(x^2-2a(t)x+(a(t)^2-4)), 
\qquad T_Y\cong \langle -8\rangle\oplus U(2)^2.
$$
This K3 surface has a Nikulin involution
defined by translation by the section
$t\mapsto (x(t),y(t))=(0,0)$ and the quotient surface is $X$.
For general $X$, the bad fibers of $Y$ are
$8$ fibers of type $I_2$ over the same points in $\PP^1$ where
$X$ has fibers of type $I_1$ and at infinity $Y$ has a fiber of type $I_8$.

\ts
As we observed in section \ref{elliptic},
translation by the section of order two defines a Nikulin involution.

Let $\hat{a}(s):=s^4a(s^{-1})$, it is a polynomial of degree at most four
and $\hat{a}(0)\neq 0$.
Then on $\PP^1-\{0\}$, with coordinate $s=t^{-1}$, the Weierstrass model is
$$
v^2=u(u^2+\hat{a}(s)u+s^8),\qquad \Delta=s^{16}(\hat{a}(s)^2-4s^8),
\qquad u=s^4x,\;v=s^6y,
$$
where $\Delta$ is the discriminant.
The fiber over $s=0$ is a stable (nodal) curve, so the corresponding fiber
$X_\infty$ is of type $I_m$ where $m$ is the order of vanishing of
the discriminant in $s=0$ (equivalently, it is the order of the pole of the
$j$-invariant in $s=0$).
Thus $X_\infty$ is an $I_{16}$ fiber.
As the section of order two specializes
to the singular point $(u,v,s)=(0,0,0)$,
after blow up it will not meet the component of the
fiber which meets the zero section.

The group structure of the elliptic fibration
induces a Lie group structure on the smooth part of the $I_{16}$
fiber.
Taking out the 16 singular points in this fiber, we
get the group $\CC^*\times \ZZ/16\ZZ$. The zero section meets the component
$C_0$, where
$$
C_n\,:=\,\PP^1\times \{\bar{n}\}\;\hookrightarrow X_\infty,
$$
and the section of order two must meet $C_8$.
Translation by the section of order two induces the permutation
$C_n\;\mapsto\; C_{n+8}$
of the $16$ components of the fiber.
The classes of the components $C_n$,
with $n=-2,\ldots,4$, generate a lattice of type $A_7(-1)$
which together with the zero section gives an $E_8(-1)$.
The Nikulin involution maps this $E_8(-1)$
to the one whose components are the
$C_n$, $n=6,\ldots,12$, and the section of order two.
Thus the Nikulin involution permutes two perpendicular copies of $E_8(-1)$
and hence it is a Morrison-Nikulin involution.

The bad fibers over $\PP^1-\{\infty\}$ correspond to the zeroes
of $\Delta=a^2(t)-4$. For general $a$, $\Delta$ has eight simple zeroes
and the fibers are nodal, so we have
eight fibers of type $I_1$ in $\PP^1-\{\infty\}$.

By considering the points on $\PP^1$ where there are bad fibers it is
not hard to see that we do get a three dimensional family of elliptic
K3 surfaces with a Morrison-Nikulin involution. Thus the general member
of this three dimensional family has a N\'eron Severi group $S$ of rank $17$.

As we constructed a unimodular sublattice $E_8(-1)^2\subset S$,
we get $S\cong <-d>\oplus E_8(-1)^2$ and $d\;(>0)$ is the discriminant of $S$.
The Shioda-Tate formula (cf.\ e.g.\ \cite[Corollary 1.7]{shioda})
gives that $d=16/n^2$ where $n$ is the order of the group of
torsion sections. As $n$ is a multiple of $2$ and $d$ must be even it
follows that $d=4$. As the embedding of $NS(X)$ into $U^3\oplus E_8(-1)^2$
is unique up to isometry it is easy to determine $T_X=NS(X)^\perp$.
Finally $T_Y\cong T_X(2)$ by the results of \cite{morrison}.

The Weierstrass model of the
quotient elliptic fibration $Y$ can be computed with the
standard formula cf.\ \cite[p.79]{ST}, the bad fibers can be found from the discriminant
$\Delta=-4(a^2-4)^2$
(and $j$-invariant). Alternatively,
fixed points of the involution on $X$
are the nodes in the $I_1$-fibers.
Since these are blown up, we get
$8$ fibers of type $I_2$ over the same points in $\PP^1$ where
$X$ has fibers of type $I_1$.
At infinity $Y$ has a fiber of type $I_8$ because the involution
on $X$ permutes of the 16 components of the $I_{16}$-fiber
($C_n\leftrightarrow C_{n+8}$).
\qed

\subsection{Remark} The Weierstrass model we used to define
$X$, $y^2=x(x^2+a(t)x+1)$,
exhibits $X$ as the minimal model of the double cover of
$\PP^1\times \PP^1$, with affine coordinates $x$ and $t$. The branch
curve consists of the
the lines $x=0$, $x=\infty$ and the curve of bidegree $(2,4)$ defined
by $x^2+a(t)x+1=0$. Special examples of such double covers are studied
in section V.23 of \cite{BPV}. In particular, on p.185 the 16-gon appears
with the two sections attached and the $E_8$'s are pointed out in the text.
Note however that our involution is {\em not} among those studied there.

\



\begin{thebibliography}{AMRT}


\bibitem[B]{barth} W.\ Barth,
{\it Even sets of eight rational curves on a $K3$-surface},
in: Complex geometry (G\"ottingen, 2000) 1--25,
Springer, Berlin, 2002.

\bibitem[BPV]{BPV} W.\ Barth, C.\ Peters, A.\ Van de Ven,
Compact complex surfaces.
Springer-Verlag, Berlin, 1984.


\bibitem[Ca]{Catanese} F.\ Catanese,
{\it Babbage's conjecture, contact of surfaces,
symmetric determinantal varieties and applications},
Invent.\ Math.\ {\bf 63} (1981) 433--465.

\bibitem[Co]{conway}
J.\ Conway, {\it Atlas of finite groups}, Oxford \ Clarendon\ Press, 1985.

\bibitem[Deg]{degtyarev} A.\ Degtyarev,
{\it On deformations of singular plane sextics},
eprint math.AG/0511379.

\bibitem[Dem]{demazure} M.\ Demazure,
{\it Surfaces de del Pezzo I-V}, S\'eminaire sur le singularit\'es des surfaces,
Cent.\ Math.\ Ac.\ Polytech, Paliseau,
LNM {\bf 777}, 21--69, Springer, 1980.

\bibitem[Do]{dolgachev} I.\ Dolgachev,
{Mirror symmetry for lattice polarized K3-surfaces},
J.\ Math.\ Sciences {\bf 81} (1996) 2599--2630.

\bibitem[DoO]{dolgort} I.\ Dolgachev, D.\ Ortland,
{\it Point sets in projective spaces and theta functions},
Ast\'erisque {\bf 165} (1988).

\bibitem[GL]{galluzzi} F.\ Galluzzi, G.\ Lombardo,
{\it Correspondences between $K3$ surfaces},
with an appendix by Igor Dolgachev,
Michigan Math.\ J.\ {\bf 52} (2004) 267--277.

\bibitem[GS]{GS} A.\ Garbagnati, A.\ Sarti,
{\it Symplectic automorphism of prime order on $K3$ surfaces},
in preparation.


\bibitem[I]{Inose}
H.\ Inose,
{\it On certain Kummer surfaces which can be realized as
non-singular quartic surfaces in $P\sp{3}$},
 J.\ Fac.\ Sci.\ Univ.\ Tokyo Sect. IA Math.\ {\bf 23} (1976) 545--560.





\bibitem[Mi]{miranda}
R.\ Miranda,
\emph{The basic theory of elliptic surfaces},
ETS Editrice  Pisa 1989.

\bibitem[Mo]{morrison}
D.R.\ Morrison,
\emph{On K3 surfaces with large Picard number},
Invent.\ Math.\ {\bf 75} (1986) 105--121.

\bibitem[Mu]{mukai} S.\ Mukai,
{\it On the moduli space of bundles on $K3$ surfaces. I},
in: Vector bundles on algebraic varieties (Bombay, 1984),
Tata Inst.\ Fund.\ Res.\ Stud.\ Math., 11, 1987, 341--413.

\bibitem[Ni1]{Nikulin} V.\ V.\ Nikulin,
{\it Finite groups of automorphisms of K\"ahlerian $K3$ surfaces},
(Russian) Trudy Moskov. Mat. Obshch. 38 (1979), 75--137, translated as:
{\it Finite automorphism groups of K\"ahler $K3$ surfaces},
Trans.\ Moscow Math.\ Soc.\ {\bf 38} (1980) 71--135.

\bibitem[Ni2]{Nikulin bilinear} V.\ V.\ Nikulin,
{\it Integral symmetric bilinear forms and some applications},
Izv.\ Math.\ Nauk SSSR {\bf 43} (1979) 111--177,
Math.\ USSR Izvestija {\bf 14} (1980) 103--167.

\bibitem[Ni3]{nikulin cor} V.\ V.\ Nikulin,
{\it On correspondences between surfaces of K3 type} (Russian)
Izv.\ Akad.\ Nauk SSSR Ser.\ Mat.\ {\bf 51} (1987) 402--411, 448;
translation in  Math.\ USSR-Izv.\ {\bf 30} (1988) 375--383.


\bibitem[P]{persson}
U.\ Persson,
\emph{Double Sextics and singular K3 surfaces},
Algebraic Geometry, Sitges (Barcelona), 1983, 262--328,
LNM 1124, Springer, Berlin, 1985.

\bibitem[R]{Reiner} I.\ Reiner,
{\it Integral representations of cyclic groups of prime order}.
Proc.\ Amer.\ Math.\ Soc.\ {\bf 8} (1957) 142--146.

\bibitem[SD]{saintdonat}
B.\ Saint-Donat, {\it Projective Models of K3 surfaces},
Amer.\ J.\ of  Math.\ {\bf 96} (1974) 602--639.

\bibitem[Se]{Serre} J.P.\ Serre, A course in Arithmetic.
GTM 7, Springer,  1973.

\bibitem[Shim]{Shimada}  I.\ Shimada, {\it On elliptic K3 surfaces},
eprint math.AG/0505140.

\bibitem[Shio]{shioda} T.\ Shioda, {\it On elliptic modular surfaces}.
J.\ Math.\ Soc.\ Japan {\bf 24} (1972) 20--59.

\bibitem[ST]{ST} J.\ Silverman, J.\ Tate,
\emph{Rational points on elliptic curves},
Undergraduate Texts in Mathematics. Springer-Verlag, New York, 1992. 


\end{thebibliography}
\end{document}